\documentclass[leqno]{article}
\usepackage{beamerarticle} 
\usepackage{geometry}      \usepackage{amsthm}        

\usepackage[T2A]{fontenc}
\usepackage[utf8]{inputenc}
\usepackage[english,russian]{babel}

\ifdefined\isenglish\else\def\isenglish{1}\fi

\if\isenglish1
  \AtBeginDocument{\selectlanguage{english}}
\else
  \AtBeginDocument{\selectlanguage{russian}}
\fi

\ifdefined\whichmode\else\def\whichmode{1}\fi

\newcommand{\mytitle}{Gromov--Hausdorff distance and stability of dynamical systems}

\newcommand{\myabstract}{
We introduce the notion of an $F$-space --- a set equipped with a family of generalized pseudometrics and marked points, and construct the categorical Gromov--Hausdorff distance for $F$-spaces.
Based on this, we propose a new definition of the distance between dynamical systems, understood in a broad sense as parameterized families of maps.
The main focus is on local dynamics: we prove the preservation of Lyapunov stability and asymptotic stability under limit transitions with respect to the new distance (in the latter case, assuming a common radius of attraction).
It is established that the class of stable systems is a closed and nowhere dense subset in the space of compact dynamical systems.
Next, we investigate the properties of the Hausdorff map, which induces an $F$-space structure on the family of subsets.
In the final part, we introduce another modified version of the Gromov--Hausdorff distance for dynamical systems and use it to calculate the exact distance between torus translations for a non-resonant vector.

Bibliography: 16 references.
}

\newcommand{\myname}{Emin~Abdullaev}

\newcommand{\mydate}{\ifnum\day<10 0\fi\the\day.\ifnum\month<10 0\fi\the\month.\the\year}

\usepackage{mathtools} \usepackage{amsfonts,amssymb,mathrsfs,amscd}
\usepackage{stmaryrd}

\usepackage{graphicx}
\usepackage{verbatim}
\usepackage{url}
\usepackage{alltt}
\usepackage{indentfirst}
\usepackage{qrcode}

\theoremstyle{plain}

\newtheorem{thm}{Theorem}[section]
\newtheorem{prop}[thm]{Proposition}
\newtheorem{cor}[thm]{Corollary}
\newtheorem{ass}[thm]{Assertion}
\newtheorem{lem}[thm]{Lemma}

\theoremstyle{definition}

\newtheorem{examp}[thm]{Example}

\newtheorem{dfn}[thm]{Definition}
\newtheorem{rk}[thm]{Remark}

\newtheorem*{agree}{Convention}
\newtheorem*{notation}{Notation}

\newcommand{\dl}{\delta}

\newcommand{\e}{\varepsilon}

\newcommand{\cGH}{\mathcal{GH}}
\newcommand{\cDM}{\mathcal{DM}}
\newcommand{\cDGH}{\mathcal{DGH}}
\newcommand{\cH}{\mathcal{H}}
\newcommand{\cR}{\mathcal{R}}
\newcommand{\cM}{\mathcal{M}}

\newcommand{\cC}{\mathcal{C}}
\newcommand{\cK}{\mathcal{K}}
\newcommand{\cF}{\mathcal{F}}
\newcommand{\cI}{\mathcal{I}}
\newcommand{\cP}{\mathcal{P}}
\newcommand{\cS}{\mathcal{S}}

\newcommand{\fB}{\mathfrak{B}}
\newcommand{\fU}{\mathfrak{U}}

\newcommand{\N}{\mathbb{N}}
\newcommand{\Q}{\mathbb{Q}}
\newcommand{\R}{\mathbb{R}}
\newcommand{\Z}{\mathbb{Z}}
\newcommand{\bU}{\mathbb{U}}

\newcommand{\diam}{{\operatorname{diam}}}
\newcommand{\Rad}{{\operatorname{Rad}}}
\newcommand{\dis}{{\operatorname{dis}}}
\newcommand{\codis}{{\operatorname{codis}}}

\newcommand{\Ob}{{\operatorname{Ob}}}
\newcommand{\Mor}{{\operatorname{Mor}}}
\newcommand{\id}{{\operatorname{id}}}

\newcommand{\Int}{{\operatorname{Int}}}

\renewcommand{\:}{\colon}
\newcommand{\rom}[1]{{\em #1}}

\makeatletter
\@ifundefined{RBibitem}{
\def\RBibitem#1{\normalfont\bibitem{#1}\ignorespaces}

\def\by{\unskip\normalfont~}
  \def\paper{\unskip\normalfont. }
  \def\jour{\unskip\normalfont, \itshape }
  \def\publaddr{\unskip\normalfont, }
  \def\publ{\unskip\normalfont, }
  \def\pages{\unskip\normalfont, pp.~}

\def\inbook{\unskip\normalfont\ \itshape }
  \def\eds{\unskip\normalfont, ed. by }
  \def\book{\unskip\normalfont. \itshape }
  \def\issue{\unskip\normalfont, no.~}

\def\yr{\unskip, \normalfont}
  \def\vol{\unskip, \normalfont\bfseries}

  \def\bysame{\by\hbox to3em{\hrulefill}\thinspace\kern\z@}
  \def\preprint{\unskip\normalfont, preprint }

  \def\finalinfo{\unskip\normalfont, }
  \def\@@finalinfo{\finalinfo}

\def\morerref{\unskip\spacefactor=1000 \newblock}
  
  \def\transl{\unskip\spacefactor=1000 \newblock Engl. transl.: }
}{}
\makeatother

\usepackage{hyperref}      
\begin{document}
\title{\mytitle}
\author{\myname}
\date{\today}

\maketitle

\begin{abstract}
\myabstract
\end{abstract}

\ifnum\whichmode=1
\tableofcontents
\fi

\section*{Introduction}

The concept of distance between metric spaces, introduced by F.~Hausdorff \cite{Hausdorff} and generalized in the works of D.~Edwards \cite{Edwards} and M.~Gromov \cite{Gromov}, is one of the fundamental tools of modern metric geometry. This theory was further developed in studies considering additional structures on metric spaces, as well as investigating the properties of the Gromov--Hausdorff distance itself (see, for example, \cite{BurBurIva, IvaNikTuz, NP01}). The concept of the continuous Gromov--Hausdorff distance received its modern treatment in~\cite{LimMemoly}. An important milestone in this area was the research of S.~A.~Bogaty and A.~A.~Tuzhilin \cite{BogatyTuzhilin} devoted to the study of its fundamental properties.

The generalization of this approach to the cases of spaces endowed with a dynamical structure and parameterized families of maps laid the foundation for the Gromov--Hausdorff theory for dynamical systems (see a detailed survey in \cite{LeeMorales}).

To study dynamical systems, we needed a structure capable of comparing exactly how the mutual arrangement of all points in the space evolves over time.
Instead of comparing the trajectories themselves, we fix each moment in time and observe how the points shift relative to each other, evaluating the distances between their images.
For a rigorous mathematical description, we introduce the concept of an $F$-space --- a set equipped with a nonempty family of generalized pseudometrics, as well as an indexed set of marked points.
By considering an arbitrary dynamical system as an $F$-space, where the family of pseudometrics is generated by translations along trajectories, we define the Gromov--Hausdorff distance between dynamical systems. Unlike in~\cite{Fukaya1986}, we do not assume the maps to be isometric. The distance we introduce satisfies the triangle inequality, unlike in~\cite{ArbietoMorales2017, Chung2020}.

In the existing literature, the Gromov--Hausdorff distance for dynamical systems is mainly used to analyze their global behavior (for example, the stability of attractors).
However, questions of local dynamics, such as the preservation of classical Lyapunov stability and asymptotic stability under limit transitions, remain poorly studied.
For a rigorous analysis of such characteristics, it is natural to consider pointed dynamical systems, and to preserve the topological structure of orbits --- to use exactly the continuous version of the Gromov--Hausdorff distance, since the usual Gromov--Hausdorff distance can ``tear'' orbits, which may cause local properties like Lyapunov stability to break down under a limit transition.
In this regard, a fundamental question arises: are classical stability properties preserved under convergence with respect to the continuous Gromov--Hausdorff distance, and what is the structure of the subspace of all stable systems?

The present work is devoted precisely to this question.
In Section 1, we introduce the concept of an $F$-space and define the categorical Gromov--Hausdorff distance between them (the word ``categorical'' merely means that the distance is calculated by taking the infimum over a certain allowed class of maps).
Next, we define the categorical (in particular, the usual and continuous) Gromov--Hausdorff distance for (pointed) dynamical systems.
In Section 2, it is proved that Lyapunov stability and asymptotic stability (in the latter case, only when a common radius of attraction exists) are preserved under convergence of systems with respect to the pointed continuous Gromov--Hausdorff distance.
It is established that the class of stable systems is closed and nowhere dense, while the subspace of asymptotically stable systems forms a nowhere dense $F_{\sigma}$-set.
In Section 3, we investigate the natural $F$-space structure on the space of all nonempty subsets with the Hausdorff distance, as well as the dynamical system structure on the Hausdorff space of closed bounded subsets, and prove the 1-Lipschitz property of the induced Hausdorff map in the general case for $F$-spaces and for continuous compact dynamical systems.
In Section 4, we return to the question of the shape of orbits. Discovering that the first (primary) version of our distance does not distinguish between some qualitatively different dynamics (e.g., ergodic torus translations), we introduce an alternative, modified version of the Gromov--Hausdorff distance for dynamical systems, which allows one to take into account not only the change in distances between points over time but also the geometry of the trajectories themselves. For this new distance, a sufficient condition for attaining the upper bound is obtained, which we illustrate using the example of torus translations.

\subsection*{Acknowledgments}
\begin{sloppypar}
The author expresses gratitude to his scientific advisor, Doctor of Physical and Mathematical Sciences, Professor \mbox{A.~A.~Tuzhilin}, and Doctor of Physical and Mathematical Sciences, Professor \mbox{A.~O.~Ivanov}, for formulating the problem and their constant attention to the work.
\end{sloppypar}

\section{Main definitions}

\subsection{Gromov--Hausdorff distance for metric spaces}

Let $\cGH$ be the class of all nonempty metric spaces, and $\cM \subset \cGH$ be the class of all nonempty metric compact spaces. For arbitrary $X, Y \in \cGH$ let $d_X$ and $d_Y$ denote their metrics.

Let $R \subset X \times Y$ be an arbitrary nonempty subset.
Define the \emph{distortion}
$$\dis(R) = \sup_{\substack{(x, y) \in R \\ (\tilde{x}, \tilde{y}) \in R}} \bigl| d_{X}(x, \tilde{x}) - d_{Y}(y, \tilde{y}) \bigr|.$$
In particular, for an arbitrary map $f \: X \to Y$ we define its distortion as the distortion of its graph:
$$\dis(f) = \sup_{x, \tilde{x} \in X} \Bigl| d_{X}(x, \tilde{x}) - d_{Y}\bigl(f(x), f(\tilde{x})\bigr)\Bigr|.$$

For a pair of maps $f_1 \: X \to Y, \ f_2 \: Y \to X$ we define their \emph{codistortion}
$$\codis(f_1, f_2) = \sup_{\substack{x \in X \\ y \in Y}} \Bigl| d_{X}\bigl(x, f_2(y)\bigr) - d_{Y}\bigl(y, f_1(x)\bigr)\Bigr|.$$

A \emph{correspondence} between $X$ and $Y$ is a subset $R \subset X \times Y$ whose projections onto $X$ and $Y$ are surjective. The set of all correspondences between $X$ and $Y$ will be denoted by $\cR(X, Y)$.
The Gromov--Hausdorff distance can be defined as
$$d_{GH}(X, Y) = \frac{1}{2} \inf_{R \in \cR(X, Y)} \dis(R).$$
In what follows, when discussing the properties of the Gromov--Hausdorff distance, we follow \cite{BurBurIva}, \cite{BogatyTuzhilin}.
The following proposition provides an equivalent definition of the Gromov--Hausdorff distance.

\begin{prop}
For any $X, Y \in \cGH$ it holds that
$$d_{GH}(X, Y) = \frac{1}{2} \inf_{\substack{f_1 \: X \to Y \\ f_2 \: Y \to X}} \max \bigl\{ \dis(f_1), \dis(f_2), \codis(f_1, f_2) \bigr\},$$
where the infimum is taken over all maps $f_1, f_2$.
\end{prop}

\subsection{Categorical Gromov--Hausdorff distance for $F$-spaces}

We will call a \emph{generalized pseudometric on a set} $X$ any mapping
$$d \: X \times X \to [0, \infty],$$ satisfying
$d(x, y) = d(y, x), \ d(x, x) = 0$ and $d(x, y) \le d(x, z) + d(z, y)$ for all $x, y, z \in X$.
If $d(x, y) = 0$ implies $x = y$, then the word ``pseudometric'' is replaced by the word \emph{metric}. If $d(x, y) < \infty$ for all $x, y \in X$, then the word ``generalized'' is omitted.

Let some index sets $I \ne \emptyset$ and $J$ be given.
We will consider these index sets fixed unless stated otherwise.
\begin{dfn}
We call an $F$-\emph{space} a nonempty set $X$ equipped with a nonempty family of generalized pseudometrics $\{d_{X, i}\}_{i \in I}$ and an arbitrary family of pairwise distinct marked points $\{x_{j}\}_{j \in J}$.
\end{dfn}
\begin{examp}
Let a family of seminorms $\{p_{i}\}_{i \in I}$ be given on a vector space $X$.
Then $X$ is an $F$-space (for example, one can assume $J = \emptyset$, or take $0$ as the marked point).
\end{examp}
Denote by $\cF_{I, J}$ the category of all $F$-spaces, where morphisms are all maps preserving marked points, that is, mapping each marked point $x_j \in X$ to the corresponding marked point $y_j \in Y$ for all $j \in J$.

In other words,
\begin{multline*}
\Mor_{\cF}\bigl((X, \{d_{X, i}\}_{i \in I}, \{x_{j}\}_{j \in J}), (Y, \{d_{Y, i}\}_{i \in I}, \{y_{j}\}_{j \in J})\bigr) = \\
= \bigl\{ f\: X \to Y \ \big|\  \forall j \in J \ f(x_j) = y_j\bigr\}.
\end{multline*}
If it is clear from the context which sets $I, J$ are considered, we will write $\cF$ instead of $\cF_{I, J}$.

Let us generalize the concept of distortion. Here and below, any difference $\bigl| \alpha - \beta \bigr|$ is understood in the extended sense: if exactly one of the values $\alpha$ or $\beta$ is infinite, then the modulus of the difference is set to $\infty$, and if both values are infinite, their difference is considered equal to $0$.

\begin{dfn}
Let $(X, d_{X})$ and $(Y, d_{Y})$ be generalized pseudometric spaces, and let $R \subset X \times Y$ be an arbitrary nonempty subset. Then the \emph{distortion} of $R$ is defined as
$$\dis_{d_{X}, d_{Y}}(R) = \sup_{(x, y), (\tilde{x}, \tilde{y}) \in R} \bigl| d_{X}(x, \tilde{x}) - d_{Y}(y, \tilde{y})\bigr|.$$
\end{dfn}
\begin{notation}
Let $\bigl(X, \{d_{X, i}\}_{i \in I}, \{x_{j}\}_{j \in J}\bigr)$, $\bigl(Y, \{d_{Y, i}\}_{i \in I}, \{y_{j}\}_{j \in J}\bigr)$ be $F$-spaces and $R \subset X \times Y$ be an arbitrary nonempty subset. For all $i \in I$ we set
$$\dis_{i}(R) := \dis_{d_{X, i}, d_{Y, i}}(R).$$
\end{notation}

\begin{dfn}
Let $\bigl(X, \{d_{X, i}\}_{i \in I}, \{x_{j}\}_{j \in J}\bigr)$, $\bigl(Y, \{d_{Y, i}\}_{i \in I}, \{y_{j}\}_{j \in J}\bigr)$ be $F$-spaces and $R \subset X \times Y$ be an arbitrary nonempty subset. Then the \emph{distortion} of $R$ is defined as
$$
\dis^{F}(R) := \sup_{i \in I} \dis_{i}(R) =
\sup_{i \in I} \sup_{(x, y), (\tilde{x}, \tilde{y}) \in R} \bigl| d_{X, i}(x, \tilde{x}) - d_{Y, i}(y, \tilde{y})\bigr|.
$$
\end{dfn}

Similarly, we generalize the concepts of map distortion and codistortion, as well as the Gromov--Hausdorff distance.
\begin{dfn}
Let $X, Y$ be $F$-spaces, and $f_1 \: X \to Y$, $f_2 \: Y \to X$ be arbitrary maps.
Then
\begin{itemize}
\item
$\dis^{F}(f_1) = \sup_{i \in I} \sup_{x, \tilde{x} \in X} \Bigl| d_{X, i}(x, \tilde{x}) - d_{Y, i}\bigl(f_1(x), f_1(\tilde{x})\bigr) \Bigr|$;
\item
$\dis^{F}(f_2) = \sup_{i \in I} \sup_{y, \tilde{y} \in Y} \Bigl| d_{Y, i}(y, \tilde{y}) - d_{X, i}\bigl(f_2(y), f_2(\tilde{y})\bigr) \Bigr|$;
\item
$\codis^{F}(f_1, f_2) = \sup_{i \in I} \sup_{\substack{x \in X \\ y \in Y}} \Bigl| d_{X, i}\bigl(x, f_2(y)\bigr) - d_{Y, i}\bigl(y, f_1(x)\bigr)\Bigr|$.
\end{itemize}
\end{dfn}

\begin{agree}
If we consider $F$-spaces $X, Y$ and no additional structures are introduced on the sets $X$ and $Y$, we will omit the superscript $F$ in the notation for the distortion of an arbitrary nonempty subset $R \subset X \times Y$, as well as for the distortions and codistortion of maps $f_1 \: X \to Y, \ f_2 \: Y \to X$:
$$\dis(R) = \dis^{F}(R), \dis(f_1) = \dis^{F}(f_1),$$
$$\dis(f_2) = \dis^{F}(f_2), \codis(f_1, f_2) = \codis^{F}(f_1, f_2).$$
\end{agree}

The following statement immediately follows from the definition.
\begin{ass}
Let $X$, $Y$ be $F$-spaces, and $R_1 \subset R_2 \subset X \times Y$ be arbitrary nonempty subsets. Then $\dis(R_1) \le \dis(R_2)$.
\end{ass}

\begin{notation}
For any sets $X, Y$ and any maps $f_1 \: X \to Y, \ f_2 \: Y \to X$, we set
$$R_{f_1, f_2} = \Bigl\{ \bigl(x, f_1(x)\bigr) \ \big| \ x \in X \Bigr\} \cup \Bigl\{\bigl(f_2(y), y\bigr) \ \big| \ y \in Y\Bigr\}.$$
\end{notation}

\begin{ass}
If $X$, $Y$ are $F$-spaces, and $f_1 \: X \to Y, f_2 \: Y \to X$ are arbitrary maps, then
$$\dis^{F}(R_{f_1, f_2}) = \max\bigl\{ \dis^{F}(f_1), \ \dis^{F}(f_2), \ \codis^{F}(f_1, f_2) \bigr\}.$$
\end{ass}

\begin{prop}
Let $X$, $Y$, $Z$ be $F$-spaces, and $f_1 \: X \to Y, f_2 \: Y \to Z$ be arbitrary maps. Then $\dis^{F}(f_2 \circ f_1) \le \dis^{F}(f_1) + \dis^{F}(f_2)$.
\end{prop}

\begin{proof}
Since for fixed $i \in I$ and $x, \tilde{x} \in X$ it holds that
\begin{multline*}
\Bigl| d_{X, i}(x, \tilde{x}) - d_{Z, i}\big(f_2 \circ f_1 (x), f_2 \circ f_1 (\tilde{x})\big)\Bigr| \le
\Bigl| d_{X, i}(x, \tilde{x}) - d_{Y, i}\big(f_1 (x), f_1 (\tilde{x})\big)\Bigr| + \\
+ \Bigl| d_{Y, i}\big(f_1 (x), f_1 (\tilde{x})\big) - d_{Z, i}\big(f_2 \circ f_1 (x), f_2 \circ f_1 (\tilde{x})\big)\Bigr| \le \dis^{F}(f_1) + \dis^{F}(f_2),
\end{multline*}
passing to the supremum over $i \in I$ and $x, \tilde{x} \in X$, we obtain the required result.
\end{proof}

\begin{prop}
Let $X$, $Y$, $Z$ be $F$-spaces, and $f_1 \: X \to Y, \ f_2 \: Y \to X, \ f'_1 \: Y \to Z, \ f'_2 \: Z \to Y$ be arbitrary maps.
Then $\codis^{F}(f'_1 \circ f_1, f_2 \circ f'_2) \le \codis^{F}(f_1, f_2) + \codis^{F}(f'_1, f'_2)$.
\end{prop}

\begin{proof}
Since for fixed $i \in I$ and $x \in X, z \in Z$ it holds that
\begin{multline*}
\Bigl| d_{X, i}(x, f_2 \circ f'_2(z)) - d_{Z, i}(z, f'_1 \circ f_1(x)) \Bigr| \le
\Bigl| d_{X, i}(x, f_2 \circ f'_2(z)) - d_{Y, i}(f_1(x), f'_2(z)) \Bigr| + \\
+ \Bigl| d_{Y, i}(f_1(x), f'_2(z)) - d_{Z, i}(z, f'_1 \circ f_1(x)) \Bigr|
\le \codis^{F}(f_1, f_2) + \codis^{F}(f'_1, f'_2),
\end{multline*}
passing to the supremum over $x \in X$, $z \in Z$, and $i \in I$, we obtain the required result.
\end{proof}

From this point until the end of the text, we will consider arbitrary subcategories $\cK, \cK' \subset \cF$.
We will always assume that all morphism sets $\Mor_{\cK}(X, Y)$ and $\Mor_{\cK'}(X, Y)$ are nonempty.
We will also always consider these categories to be fixed, unless stated otherwise in what follows.

\begin{dfn}
Let $X, Y \in \Ob(\cK)$. Then the (\emph{categorical}) \emph{Gromov--Hausdorff distance} is defined as
$$d_{GH, \cK}^{F}(X, Y) = \frac{1}{2}
\inf_{\substack{f_1 \in \Mor_{\cK}(X, Y) \\ f_2 \in \Mor_{\cK}(Y, X)}} \max\bigl\{ \dis^{F}(f_1), \ \dis^{F}(f_2), \ \codis^{F}(f_1, f_2) \bigr\}.$$
\end{dfn}

\begin{agree}
By analogy with the distortion and codistortion, we will omit the superscript $F$ in the notation $d_{GH, \cK}^{F}(X, Y)$ if it is clear from the context that we are considering the Gromov--Hausdorff distance specifically for $F$-spaces.
\end{agree}

\begin{ass}\label{ass:metric}
The Gromov--Hausdorff distance $d^{F}_{GH, \cK}$ is a generalized pseudometric on $\Ob(\cK)$.
\end{ass}

\begin{proof}
The equality $d^{F}_{GH, \cK}(X, Y) = d^{F}_{GH, \cK}(Y, X)$ is obvious.
The equality $d^{F}_{GH, \cK}(X, X) = 0$ is explained by the fact that the distortion of the identity map is $\dis^{F}(\id_{X}) = 0$.

Let us verify the triangle inequality. Let $X, Y, Z \in \cF$. We will assume that $d^{F}_{GH, \cK}(X, Y) < \infty$ and $d^{F}_{GH, \cK}(Y, Z) < \infty$, since otherwise there is nothing to prove.
Fix $\e > 0$.
Consider morphisms $f_1 \in \Mor_{\cK}(X, Y), \ f_2 \in \Mor_{\cK}(Y, X)$ and $f'_1 \in \Mor_{\cK}(Y, Z), \ f'_2 \in \Mor_{\cK}(Z, Y)$, such that
$\dis^{F}(R_{f_1, f_2}) \le 2d^{F}_{GH, \cK}(X, Y) + \e$ and $\dis^{F}(R_{f'_1, f'_2}) \le 2d^{F}_{GH, \cK}(Y, Z) + \e$.
Then
\begin{multline*}
2d^{F}_{GH, \cK}(X, Z) \le \dis^{F}(R_{f'_1 \circ f_1, f_2 \circ f'_2}) = \\
= \max\bigl\{ \dis^{F}(f'_1 \circ f_1), \ \dis^{F}(f_2 \circ f'_2), \ \codis^{F}(f'_1 \circ f_1, f_2 \circ f'_2) \bigr\} \le \\
\le \max\bigl\{ \dis^{F}(f_1), \ \dis^{F}(f_2), \ \codis^{F}(f_1, f_2) \bigr\} + \\
\max\bigl\{ \dis^{F}(f'_1), \ \dis^{F}(f'_2), \ \codis^{F}(f'_1, f'_2) \bigr\} = \\
= \dis^{F}(R_{f_1, f_2}) + \dis^{F}(R_{f'_1, f'_2}) \le 2d^{F}_{GH, \cK}(X, Y) + 2d^{F}_{GH, \cK}(Y, Z) + 2\e.
\end{multline*}
Since $\e > 0$ is arbitrary, we obtain
$d^{F}_{GH, \cK}(X, Z) \le d^{F}_{GH, \cK}(X, Y) + d^{F}_{GH, \cK}(Y, Z)$.
\end{proof}

\begin{rk}
Note that in the case $\cK = \cF$, for any $F$-spaces $X, Y$ it holds that
$$d^{F}_{GH, \cK}(X, Y) = d^{F}_{GH, \cF}(X, Y) = \frac{1}{2} \inf_{R \in \cR(X, Y)} \dis^{F}(R).$$
\end{rk}

The following statement immediately follows from the definition.

\begin{prop}
Suppose that $\cK' \subset \cK$. Then for any $X, Y \in \Ob(\cK')$ it holds that
$$d^{F}_{GH, \cK}(X, Y) \le d^{F}_{GH, \cK'}(X, Y).$$
\end{prop}

\begin{dfn}
Let $X$ be an $F$-space. For all $i \in I$, we define the \emph{diameters}
$$\diam_{i}(X) = \sup_{x, \tilde{x} \in X} d_{X, i}(x, \tilde{x}), \ \diam(X) = \sup_{i \in I}\diam_{i}(X).$$
If $J \ne \emptyset$, then for all $i \in I, j \in J$ we define the \emph{radii}
$$\Rad_{i,j}(X) = \sup_{x \in X}d_{X, i}(x, x_j).$$
An $F$-space $X$ is called \emph{bounded} if $\diam(X) < \infty$.
\end{dfn}

\begin{prop}\label{prop:estimation}
\begin{enumerate}
\item
For any $X, Y \in \Ob(\cK)$, the following estimate holds
$$\sup_{i \in I} \ \bigl| \diam_{i}(X) - \diam_{i}(Y) \bigr| \le 2d^{F}_{GH, \cK}(X, Y) \le
\max\bigl\{ \diam(X), \diam(Y) \bigr\}.$$
\item
If $J \ne \emptyset$, then for any $X, Y \in \Ob(\cK)$, the following inequality holds
$$\sup_{\substack{i \in I \\ j \in J}} \bigl| \Rad_{i, j}(X) - \Rad_{i, j}(Y) \bigr| \le
2d^{F}_{GH, \cK}(X, Y).$$
\end{enumerate}
\end{prop}

\begin{proof}
We will assume that $X$ and $Y$ are bounded (the other cases are similar).
\begin{enumerate}
\item
Let us prove the first inequality.
Fix $i \in I$ such that $\diam_{i}(X) \ne \diam_{i}(Y)$. Without loss of generality, we will assume that $\diam_i(Y) < \diam_i(X)$. Fix $\e > 0, \ \e < \diam_{i}(X) - \diam_{i}(Y)$. Let $x, \tilde{x} \in X$ be such that $d_{X, i}(x, \tilde{x}) > \diam_{i}(X) - \e$.
Consider $f_1 \: X \to Y$, an arbitrary map.
Since
$$d_{X, i}(x, \tilde{x}) > \diam_{i}(X) - \e \text{ and } d_{Y, i}\bigl(f_{1}(x), f_{1}(\tilde{x})\bigr) \le \diam_{i}(Y),$$
we have $\dis(f_1) \ge \bigl| \diam_{i}(X) - \diam_{i}(Y) \bigr| - \e$.
Letting $\e \to 0$, we obtain $\dis(f_1) \ge \bigl| \diam_{i}(X) - \diam_{i}(Y) \bigr|$.
From this it follows that $2d_{GH, \cK}(X, Y) \ge \bigl| \diam_{i}(X) - \diam_{i}(Y) \bigr|$.

Let us prove the second inequality.
Set $R = X \times Y$. Let
$$f_1 \in \Mor_{\cK}(X, Y), \ f_2 \in \Mor_{\cK}(Y, X).$$
Then
\begin{multline*}
2d_{GH, \cK}(X, Y) \le \dis(R_{f_1, f_2}) \le \dis(R)
= \sup_{i \in I} \sup_{\substack{x, \tilde{x} \in X \\ y, \tilde{y} \in Y}}
\bigl| d_{X, i}(x, \tilde{x}) - d_{Y, i}(y, \tilde{y})\bigr| = \\
= \sup_{i \in I} \max\bigl\{ \diam_{i}(X), \diam_{i}(Y) \bigr\} = \max\bigl\{ \diam(X), \diam(Y)\bigr\}.
\end{multline*}

\item
Let $\{x_{j}\}_{j \in J}$ and $\{y_{j}\}_{j \in J}$ denote the families of marked points in $X$ and $Y$, respectively.
Fix $i \in I$ and $j \in J$ such that $\Rad_{i, j}(X) \ne \Rad_{i, j}(Y)$.
Without loss of generality, we may assume that $\Rad_{i, j}(X) > \Rad_{i, j}(Y)$.
Take $\e > 0$ such that $\e < \Rad_{i, j}(X) - \Rad_{i, j}(Y)$.
Let $f_1 \: X \to Y$ be an arbitrary map.
Choose a point $x \in X$ such that $d_{X, i}(x, x_{j}) > \Rad_{i, j}(X) - \e$.
Since $d_{Y, i}\bigl(f_1(x), y_{j}\bigr) \le \Rad_{i, j}(Y)$, we have
$\dis(f_1) \ge \bigl| \Rad_{i, j}(X) - \Rad_{i, j}(Y) \bigr| - \e.$
Since this is true for any $\e > 0$ and any $i \in I, \ j \in J$, we get
$2d_{GH, \cK}(X, Y) \ge \sup_{i \in I, \ j \in J} \bigl| \Rad_{i, j}(X) - \Rad_{i, j}(Y) \bigr|$.
\end{enumerate}
\end{proof}

Note that the Gromov--Hausdorff distance for (pseudo)metric spaces is a special case of the categorical distance for $F$-spaces we introduced (in this case, for $I = \{0\}, J = \emptyset$, the category $\cF$ of all pseudometric spaces is considered as $\cK$).
Thus, the proved propositions are also valid for the Gromov--Hausdorff distance $d_{GH}$.

Let an arbitrary $F$-space $X$ and $\lambda > 0$ be given. Then $\lambda X$ denotes the $F$-space on the set $X$ with the family of generalized pseudometrics $d_{\lambda X, i}(x_1, x_2) = \lambda d_{X, i}(x_1, x_2)$.
Also, for any $\lambda, \mu > 0$ and any $F$-space $X$, we denote by $\id_{\lambda, \mu, X} \: \lambda X \to \mu X$ the identity map, $\id_{\lambda, \mu, X}(x) = x$.

\begin{prop}
\begin{enumerate}
\item
Suppose that for any $\lambda > 0$ there is a functor
$\Lambda_{\lambda} \: \cK \to \cK, \ X \mapsto \lambda X$, which acts identically on morphisms:
$$\Lambda_{\lambda} \: \Mor_{\cK}(X, Y) \to \Mor_{\cK}(\lambda X, \lambda Y), \ f \mapsto f.$$
Then for any $X, Y \in \Ob(\cK)$ and any $\lambda > 0$ it holds that
$$d_{GH, \cK}(\lambda X, \lambda Y) = \lambda d_{GH, \cK}(X, Y).$$
\item
Suppose that for all $X \in \Ob(\cK)$ and all $\lambda, \mu > 0$, the map $\id_{\lambda, \mu, X}$ is an element of $\Mor_{\cK}(\lambda X, \mu X)$.
Then for any bounded $X \in \Ob(\cK)$ and any $\lambda, \mu \ge 0$ it holds that
$$2d_{GH, \cK}(\lambda X, \mu X) = |\lambda - \mu| \ \diam(X).$$
\end{enumerate}
\end{prop}

\begin{proof}
\begin{enumerate}
\item
Let $f_1 \in \Mor_{\cK}(X, Y), \ f_2 \in \Mor_{\cK}(Y, X)$.
Then $\dis(R_{\Lambda_{\lambda}(f_1), \Lambda_{\lambda}(f_2)}) = \lambda \dis(R_{f_1, f_2})$.
From this we obtain $d_{GH, \cK}(\lambda X, \lambda Y) = \lambda d_{GH, \cK}(X, Y)$.
\item
By~\ref{prop:estimation} we have $2d_{GH, \cK}(\lambda X, \mu X) \ge \bigl|\diam(\lambda X)-\diam(\mu X)\bigr| = |\lambda - \mu| \ \diam(X)$.
Now consider $f_1 = \id_{\lambda, \mu, X}, \ f_2 = \id_{\mu, \lambda, X}$ --- identity maps.
Then $\dis(R_{f_1, f_2}) = |\lambda - \mu| \ \diam(X)$.
Consequently, $2d_{GH, \cK}(\lambda X, \mu X) = |\lambda - \mu| \ \diam(X)$.
\end{enumerate}
\end{proof}

We state the following observation, which will be useful to us later.

\begin{ass}\label{ass:map}
Let arbitrary sets $X, Y$ and a map $f \: X \to Y$ be given.
Let a generalized pseudometric $d_{Y}$ be given on the set $Y$.
Set $d_{X}(x, y) = d_{Y}\bigl(f(x), f(y)\bigr)$.
\begin{enumerate}
\item
Then $d_{X}$ is a generalized pseudometric.
\item
If $d_{Y}$ is a metric and the map $f$ is injective, then $d_{X}$ is a metric.
\end{enumerate}
\end{ass}

\begin{agree}
We will denote the Gromov--Hausdorff distance for pseudometric spaces (where we impose no restrictions on morphisms) by $d_{GH}$, exactly as for metric spaces.
\end{agree}

\begin{prop}\label{prop:EstimationByOrdinaryGH}
For any $F$-spaces
$$\bigl(X, \{d_{X, i}\}_{i \in I}, \{x_{j}\}_{j \in J}\bigr), \bigl(Y, \{d_{Y, i}\}_{i \in I}, \{y_{j}\}_{j \in J}\bigr)$$
the following inequality holds
$$d^{F}_{GH, \cF}(X, Y) \ge \sup_{i \in I} d_{GH}\bigl((X, d_{X, i}), (Y, d_{Y, i})\bigr).$$
\end{prop}
\begin{proof}
Indeed, for any $\emptyset \ne R \subset X \times Y$ and for all $i \in I$, we have $\dis^{F}(R) \ge \dis_{i}(R)$, which immediately implies the required inequality.
\end{proof}

Note that the inequality in Proposition~\ref{prop:EstimationByOrdinaryGH} cannot be replaced by equality, as illustrated by the following example.

\begin{examp}
Let $I$ be the set of all bijections $f \: \R \to \R$, let $J = \emptyset$.
Consider the $F$-space $X = \R$ with the family of metrics $d_{X, i}(x, y) = |x - y|$, $i \in I$. As the space $Y$ we consider the set $\R$, on which the family of metrics $\{d_{Y, i}\}_{i \in I}$ is defined as follows: for each bijection $f \in I$, $f \: \R \to \R$ we set $d_{Y, f}(x, y) = \bigl|f(x) - f(y)\bigr|$.
Note that for all $f \in I$ we have $\dis_{d_{X, f}, d_{Y, f}}(f) = 0$, which implies $$\sup_{i \in I} d_{GH}\bigl((X, d_{X, i}), (Y, d_{Y, i})\bigr) = 0.$$
Now we compute the Gromov--Hausdorff distance for the $F$-spaces $X, Y$.
Fix an arbitrary correspondence $R \subset X \times Y$ and a number $0 < M < \infty$. Consider any $y_1, y_2 \in Y$, for which $y_1 \ne y_2$, and $x_1, x_2 \in X$ such that $(x_1, y_1), (x_2, y_2) \in R$.
Consider a bijection $f \in I$ such that $\bigl|f(y_1) - f(y_2)\bigr| > M + |x_1 - x_2|$.
Then $\dis^{F}(R) \ge \dis_{d_{X, f}, d_{Y, f}}(R) \ge |d_{X, f}(x_1, x_2) - d_{Y, f}(y_1, y_2)| > M$. Since this holds for any correspondence $R$ and any $M > 0$, we get
$$d^{F}_{GH, \cF}(X, Y) = \infty.$$
\end{examp}

\subsection{Dynamical systems}

\begin{dfn}
A \emph{pointed} (\emph{generalized}) (\emph{pseudo})\emph{metric space} is a (generalized) (pseudo)metric space $(X, d_X)$ with a chosen base point $x_0$. We will denote it by $(X, d_X, x_0)$.
\end{dfn}

From now on, we fix a set $G$ with a distinguished point $e \in G$.

\begin{notation}
For a generalized pseudometric space $(X, d_{X})$ and any point $x \in X$ we denote by $B_{X}(x, r)$ the open ball centered at $x$ of radius $r$.
\end{notation}

\begin{dfn}
We define a \emph{dynamical system} $(X, d_X, \phi_{X})$ as a non-empty set $X$ equipped with a generalized pseudometric $d_X$ and a map $\phi_{X} \: G \times X \to X$, which for each $g \in G$ is written as $x \mapsto gx$ and satisfies the following condition:
for every $x \in X$, $ex = x$.
\end{dfn}

Note that in our definition we do not require $G$ to possess an algebraic structure (a group or a semigroup). Thus, by a dynamical system in this paper, we mean a parameterized family of maps in the broadest sense.

\begin{examp}
Let $(X, d_{X})$ be a closed smooth Riemannian manifold, and $v$ be a smooth vector field.
By the classical existence and uniqueness theorem, such a field generates a global one-parameter phase flow, which is a family of diffeomorphisms $\{g^t\}_{t \in \R}$, where $g^t \: X \to X$.
Let $G = \R$ with the distinguished point $e = 0$.
Then the dynamical system is defined as $(X, d_X, \phi_{X})$, where the action $\phi_{X} \: \R \times X \to X$ is given by the phase flow
$$\phi_{X}(t, x) = g^{t}(x).$$
The condition $ex = x$ is satisfied automatically, since $g^{0}(x) = x$.
\end{examp}

\begin{examp}
Let $(X, d_{X})$ be an arbitrary metric space endowed with a semigroup of maps $S^{t} \: X \to X, \ t \ge 0$, that is, $S^{0} = \id_{X}$ and $S^{t_1 + t_2} = S^{t_1} \circ S^{t_2}$ for all $t_1, t_2 \ge 0$.
In this case, the structure of the dynamical system $(X, d_{X}, \phi_{X})$ is defined as follows.
Let $G = [0, \infty)$ with the distinguished point $e = 0$.
The action $\phi_{X} \: [0, \infty) \to X$ is defined as $\phi_{X}(t, x) = S^{t}(x)$.
The condition $ex = S^{0}(x) = x$ is satisfied by the definition of the semigroup identity.
\end{examp}

\begin{examp}
Let $G = \R^{n}$ with the distinguished point $e = 0$.
Let $(X, d_{X})$ be a closed smooth Riemannian manifold endowed with a set of smooth vector fields $v_1, \dots, v_n$, in involution (that is, all their commutators vanish: $[v_i, v_j] = 0$ for $i \ne j$).
It is known (see, e.g., Theorem 9.44 in~\cite{Lee2013}) that in this case the one-parameter phase flows generated by them, $\{g^{t_{1}}_{1}\}_{t_{1} \in \R}, \dots, \{g^{t_{n}}_{n}\}_{t_{n} \in \R}$ commute pairwise: $$g_{i}^{t_{i}} \circ g_{j}^{t_{j}} = g_{j}^{t_{j}}\circ g_{i}^{t_{i}}.$$
This allows us to define the action as the composition of the phase flows
$$\phi_{X} \: \R^{n} \times X \to X, \ \phi_{X}\bigl((t_1, \dots, t_n), x \bigr) = (g^{t_1} \circ \dots \circ g^{t_n})(x).$$
Thus, $(X, d_{X}, \phi_{X})$ is a dynamical system in the sense of our definition.

Note that the trajectory $\phi_{X}(t, x)$ is the unique solution to the system of differential equations
$$\frac{\partial \phi_{X}(t, x)}{\partial t_k} = v_{k}\bigl(\phi_{X}(t, x)\bigr), \ k = 1, \dots, n$$ with the initial condition $\phi_{X}(0, x) = x$,
and the assumption $[v_i, v_j] = 0$ serves as the compatibility condition for the system.

\end{examp}

\begin{dfn}
We define a \emph{pointed dynamical system} as a dynamical system with a distinguished point $x_0$.
\end{dfn}

\begin{dfn}
The point $x_0$ of a pointed dynamical system $(X, x_0)$ is called an \emph{equilibrium point} if for all $g \in G$ we have $gx_0 = x_0$.
\end{dfn}

For any dynamical system $(X, d_{X}, \phi_{X})$, according to~\ref{ass:map}, we can introduce the following set of generalized pseudometrics $\{d_{X, g}\}_{g \in G}$:
$$d_{X, g}(x_1, x_2) = d_{X}(gx_1, gx_2).$$

Thus, each dynamical system $X$ defines an associated $F$-space with $I = G$, $J = \emptyset$.

We will consider a pointed dynamical system as an $F$-space with $I = G$, $J = \{0\}$.

\begin{agree}
We will say that a certain property of a (pointed) (generalized) (pseudo)metric space holds for a (pointed) dynamical system if it holds for the corresponding (pointed) (generalized) (pseudo)metric space.
\end{agree}

\begin{rk}
Note that the diameter of a dynamical system $X$ in the sense of an $F$-space is
$$\diam(X) = \sup_{g \in G}\diam_{g}(X) = \diam_{e}(X)$$
and is equal to the diameter in the sense of a generalized pseudometric space.
\end{rk}

\begin{rk}
Note that the radius $\Rad_{e, 0}(X, x_0)$ of a pointed dynamical system $(X, x_0)$ in the sense of an $F$-space is the radius of the generalized pseudometric space $X$; we will denote it for brevity by $\Rad(X, x_0)$.
\end{rk}

\subsection{Gromov--Hausdorff distance for dynamical systems}

Denote by $\cDGH = \cDGH^{G}$ the class of all dynamical systems, and by $\cDGH_{p} = \cDGH^{G}_{p}$ the class of all pointed dynamical systems.

According to the construction described at the end of the previous subsection, each dynamical system (respectively, pointed dynamical system) is naturally endowed with the structure of an $F$-space.
Denote this canonical embedding of classes by $\cI$. This embedding associates:

\begin{enumerate}
\item
to an arbitrary dynamical system $(X, d_{X}, \phi_{X})$ from the class $\cDGH$ the associated $F$-space
$$\cI(X) = \bigl(X, \{d_{X, g}\}_{g \in G}\bigr) \in \cF_{G, \emptyset};$$
\item
to an arbitrary pointed dynamical system $(X, d_{X}, x_0, \phi_{X})$ from the class $\cDGH_{p}$ the associated $F$-space
$$\cI(X) = \bigl(X, \{d_{X, g}\}_{g \in G}, x_0 \bigr) \in \cF_{G, \{0\}}.$$
\end{enumerate}

For the correctness of further definitions, we impose the following requirement on the subcategory $\cK \subset \cF$: the class of its objects $\Ob(\cK)$ must be closed under this association. Namely, we assume that
$\cI(\cDGH) \subset \Ob(\cK)$ (or $\cI(\cDGH_{p}) \subset \Ob(\cK)$ in the case of pointed systems).
This allows one to consider (pointed) dynamical systems as a special case of objects of the category $\cK$.

\begin{agree}
To avoid ambiguity, we emphasize that the subcategory $\cK \subset \cF$ is understood depending on the mathematical context as follows.
\begin{enumerate}
\item
When computing the distance between dynamical systems (from the class $\cDGH$) we assume a priori that the category $\cF$ is fixed as a subcategory in $\cF_{I, \emptyset}$, that is, we consider $F$-spaces without distinguished points ($J = \emptyset$).
\item
When passing to pointed dynamical systems (from the class $\cDGH_{p}$) the symbol $\cK$ always implies the corresponding subcategory in $\cF_{I,\{0\}}$, that is, its objects are equipped with one distinguished point ($J = \{0\}$), and morphisms are required to map one distinguished point to another.
\end{enumerate}
\end{agree}

Now we can define the Gromov--Hausdorff distance between (pointed) dynamical systems as the distance between the corresponding $F$-spaces. In other words, for any (pointed) dynamical systems $X$ and $Y$ we set
$$d^{G}_{GH, \cK}(X, Y) := d^{F}_{GH, \cK}\bigl(\cI(X), \cI(Y)\bigr).$$

\begin{notation}
Let $X, Y$ be (pointed) dynamical systems considered as $F$-spaces by means of the canonical embedding $\cI$.
To emphasize the dynamical nature of the constructions, we will add the superscript $G$ to the notation for distortion and codistortion. Namely, for arbitrary maps $f_1: X \to Y$ and $f_2: Y \to X$ their dynamical distortion and codistortion are defined as the corresponding values for the associated $F$-spaces:
\begin{itemize}
\item
$\dis^{G}(f_1) := \sup_{g \in G} \ \sup_{x, \tilde{x} \in X} \
\Bigl| d_{X}(gx, g\tilde{x}) - d_{Y}\bigl(gf_1(x), gf_1(\tilde{x})\bigr)\Bigr|$;
\item
$\dis^{G}(f_2) := \sup_{g \in G} \ \sup_{y, \tilde{y} \in Y} \
\Bigl| d_{Y}(gy, g\tilde{y}) - d_{X}\bigl(gf_2(y), gf_2(\tilde{y})\bigr)\Bigr|$;
\item
$\codis^{G}(f_1, f_2) := \sup_{g \in G} \sup_{x \in X, y \in Y} \Bigl|d_{X}\bigl(gx, gf_2(y)\bigr) - d_{Y}\bigl(gy, gf_1(x)\bigr)\Bigr|$.
\end{itemize}
\end{notation}

\begin{rk}\label{rk:AnotherDfn}
Note that for any (pointed) dynamical systems $X, Y$ we have
\begin{multline*}
d^{G}_{GH, \cK}(X, Y) = \frac{1}{2} \inf_{\substack{f_1 \in \Mor_{\cK}(\cI(X), \cI(Y)) \\ f_2 \in \Mor_{\cK}(\cI(Y), \cI(X))}} \dis^{G}(R_{f_1, f_2}) = \\ = \frac{1}{2} \inf_{\substack{f_1 \in \Mor_{\cK}(\cI(X), \cI(Y)) \\ f_2 \in \Mor_{\cK}(\cI(Y), \cI(X))}} \max\bigl\{ \dis^{G}(f_1), \ \dis^{G}(f_2), \ \codis^{G}(f_1, f_2) \bigr\}.
\end{multline*}
\end{rk}

\begin{notation}
For any dynamical systems $(X, d_{X}, \phi_{X}), \ (Y, d_{Y}, \phi_{Y})$ denote by $\cC(X, Y)$ the set of all continuous maps from $(X, d_{X})$ to $(Y, d_{Y})$.
For any pointed dynamical systems
$$(X, d_{X}, x_{0}, \phi_{X}), \ (Y, d_{Y}, y_{0}, \phi_{Y})$$
we set
$$\cC_{p}(X, Y) = \cC_{p}\bigl((X, x_0), (Y, y_0)\bigr) = \bigl\{ f \in \cC(X, Y) \ | \ f(x_0) = y_0 \bigr\}.$$
\end{notation}

\begin{dfn}
\begin{enumerate}
\item
Let $J = \emptyset$.
\begin{enumerate}
\item \label{dfn:GH}
Let $\cK = \cF$.
Then for any dynamical systems $X, Y$
we define the Gromov--Hausdorff distance
$$d^{G}_{GH}(X, Y) := d^{G}_{GH, \cF}(X, Y) = \frac{1}{2}
\inf_{\substack{f_1 \: X \to Y \\ f_2 \: Y \to X}} \dis^{G}(R_{f_1, f_2}).$$
\item \label{dfn:GHc}
Now let $\Mor_{\cK}(X, Y) = \cC(X, Y)$.
Then for any dynamical systems $X, Y$ we define the \emph{continuous Gromov--Hausdorff distance}
$$d^{G}_{GH, c}(X, Y) := d^{G}_{GH, \cK}(X, Y) = \frac{1}{2} \inf_{\substack{f_1 \in \cC(X, Y) \\ f_2 \in \cC(Y, X)}} \dis^{G}(R_{f_1, f_2}).$$
\end{enumerate}
\item
Let $J = \{0\}$.
\begin{enumerate}
\item \label{dfn:GHp}
Let $\cK = \cF$.
Then for any pointed dynamical systems $(X, x_0), (Y, y_0)$
we define the \emph{pointed Gromov--Hausdorff distance}
$$d^{G}_{GH, p}\bigl((X, x_0), (Y, y_0)\bigr) := d^{G}_{GH, \cK}(X, Y) = \frac{1}{2}
\inf_{\substack{f_1 \: X \to Y, \ f_1(x_0) = y_0 \\ f_2 \: Y \to X, \ f_2(y_0) = x_0}} \dis^{G}(R_{f_1, f_2}).$$
\item \label{dfn:GHpc}
Now let $\Mor_{\cK}(X, Y) = \cC_{p}(X, Y)$. Then for any pointed dynamical systems $(X, x_0), (Y, y_0)$ we define the \emph{pointed continuous Gromov--Hausdorff distance}
$$d^{G}_{GH, p, c}\bigl((X, x_0), (Y, y_0)\bigr) := d^{G}_{GH, \cK}(X, Y) = \frac{1}{2} \inf_{\substack{f_1 \in \cC_{p}(X, Y) \\ f_2 \in \cC_{p}(Y, X)}} \dis^{G}(R_{f_1, f_2}).$$
\end{enumerate}
\end{enumerate}
\end{dfn}

\begin{notation}
For convenience, we fix the following notation for the classes of non-empty dynamical systems with action $G$:
\begin{itemize}
\item $\cDGH$ --- the class of all dynamical systems endowed with the Gromov--Hausdorff distance $d^{G}_{GH}$, $\cDM \subset \cDGH$ --- the subclass of all compact dynamical systems with the same distance;
\item $\cDGH_{c}$ --- the class of all dynamical systems endowed with the continuous Gromov--Hausdorff distance $d^{G}_{GH,c}$, $\cDM_{c} \subset \cDGH_{c}$ --- the subclass of all compact dynamical systems with the same distance;
\item $\cDGH_{p}$ --- the class of all pointed dynamical systems endowed with the pointed Gromov--Hausdorff distance $d^{G}_{GH,p}$, $\cDM_{p} \subset \cDGH_{p}$ --- the subclass of all compact pointed dynamical systems with the same distance;
\item $\cDGH_{p,c}$ --- the class of all pointed dynamical systems endowed with the pointed continuous Gromov--Hausdorff distance $d^{G}_{GH,p,c}$, $\cDM_{p,c} \subset \cDGH_{p,c}$ --- the subclass of all compact pointed dynamical systems with the same distance.
\end{itemize}
\end{notation}

\begin{rk}
Note that any (pointed) generalized pseudometric space $X$ can be considered as a dynamical system $X^{\id}$ with the action
$$\phi^{\id}_{X} \: G \times X \to X, \ \phi^{\id}_{X}(g, x) = x.$$
Then we can introduce the Gromov--Hausdorff distance between any (pointed) generalized pseudometric spaces $X, Y$ as
$$d_{GH, \cK}(X, Y) := d^{G}_{GH, \cK}(X^{\id}, Y^{\id}).$$
The distances $d_{GH}$ and $d_{GH,c}$ defined in this way for metric spaces agree with Subsection 1.1 and~\cite{BogatyTuzhilin}.

Also note that Remark~\ref{rk:AnotherDfn} and the inequality $\dis^{G}(R) \ge \dis(R)$ imply that for any (pointed) dynamical systems $X, Y$ we have
$$d^{G}_{GH, \cK}(X, Y) \ge d_{GH, \cK}(X, Y) := d^{G}_{GH, \cK}(X^{\id}, Y^{\id}).$$
\end{rk}

\subsection{Limit of a function along a filter base}

Let us recall the concept of a filter base (base of sets) and the limit of a function along a base, for more details see~\cite{Bourbaki1958topology}.

\begin{dfn}
A \emph{filter base} (or a base of sets) on a non-empty set $X$ is a family of its subsets $\fB \ne \emptyset$, satisfying the following conditions:
\begin{enumerate}
\item $\emptyset \notin \fB$;
\item $\forall B_1, B_2 \in \fB \ \exists B_3 \in \fB \ (B_3 \subset B_1 \cap B_2)$.
\end{enumerate}
\end{dfn}

\begin{dfn}
Let a filter base $\fB$ be given on a set $X$, and $Y$ be an arbitrary topological space. A point $A \in Y$ \emph{is called the limit of the map} $f \: X \to Y$ along the base $\fB$ (denoted by $\lim_{\fB} f(x) = A$ or $f(x) \to_{\fB} A$) if for any neighborhood $U \subset Y$ of the point $A$ there exists a base element $B \in \fB$ for which $f(B) \subset U$.
\end{dfn}

\begin{rk}
If the space $Y$ is Hausdorff, then the limit of the map $f$ along the base $\fB$ is unique.
\end{rk}

Note that if $Y$ coincides with $[0, \infty]$ with the natural topology, then for a finite point $A \ge 0$ this definition can be rewritten as follows:
$$\lim_{\fB} f(x) = A \iff \forall \e > 0 \ \exists B \in \fB \ \forall x \in B \ \bigl|f(x) - A\bigr| < \e.$$

\section{Lyapunov stability}

\subsection{Lyapunov stability}

Let us generalize the definition of Lyapunov stability.
In what follows, we will assume that on the previously fixed set $G$ with the distinguished point $e$ there is also given a filter base $\fB_{G}$.

\begin{notation}
For any topological space $X$ and any of its subsets $A \subset X$ we will denote by $\Int(A)$ the interior of the set $A$.
\end{notation}

\begin{dfn}
A pointed dynamical system $(X, d_{X}, x_0)$ is called \emph{stable} (\emph{in the sense of Lyapunov}) if for any $\e > 0$ there exists $\delta > 0$ such that for all $x$ satisfying $d_{X}(x, x_0) < \delta$ it holds that for all $g \in G$ we have $d_{X}(gx, gx_0) < \e$. In other words,
$$\forall \e > 0 \ \exists \dl > 0 \ \forall x \in X \ \forall g \in G \ \bigl( d_{X}(x, x_0) < \dl \implies d_{X}(gx, gx_0) < \e \bigr). $$
\end{dfn}

\begin{dfn}
\emph{The basin of attraction} of a dynamical system $X$ with respect to a point $x_0 \in X$ is defined as the set
$$A_{X, x_0} := \Bigl\{x \in X \ \Big| \ \lim_{\fB_{G}} d_{X}(gx, gx_0) = 0 \Bigr\}.$$
\end{dfn}

It follows from the definition that we always have $x_0 \in A_{X, x_0}$.

\begin{dfn}
A pointed dynamical system $(X, d_{X}, x_0, \phi_{X})$ is called \emph{quasi-asymptotically stable} if the point $x_0$ is an interior point of the basin of attraction: $x_0 \in \Int(A_{X, x_{0}})$.
In this case, a number $\dl > 0$ is called a \emph{radius of attraction} if the point $x_0$ is contained in the set $A_{X, x_0}$ together with the open ball $B_{X}(x_0, \dl)$.
\end{dfn}

\begin{dfn}
A pointed dynamical system $X$ is called \emph{asymptotically stable} if it is both stable and quasi-asymptotically stable.
\end{dfn}

\begin{thm}\label{thm:Stability}
Let a sequence of Lyapunov stable pointed dynamical systems $\{(X_n, d_{X_n}, x^{0}_{n})\}_{n \in \N}$ converge with respect to the continuous Gromov-Hausdorff distance to a pointed dynamical system $(X, d_X, x^{0})$\/\rom:
$$d^{G}_{GH, p, c}(X_n, X) \to 0.$$
Then the pointed dynamical system $(X, d_X, x^{0})$ is Lyapunov stable.
\end{thm}
\begin{proof}
Fix $\e > 0$. Consider a dynamical system $(X_N, d_{X_N}, x^{0}_{N})$ such that $2d_{GH,p,c}\bigl((X, x^{0}), (X_N, x^{0}_{N})\bigr) < \e$.
Then there exists a continuous map
$f_1 \: X \to X_N, \ x^{0} \mapsto x^{0}_{N}$ with distortion $\dis^{G}(f_1) < 2\e$.
By the definition of stability for the system $(X_N, x^{0}_{N})$ consider $\delta > 0$ such that for all $x \in X_N$ satisfying $d_{X_N}(x, x^{0}_{N}) < \delta$ we have $d_{X_N}(gx, gx^{0}_{N}) < \e$ for all $g \in G$.
Since the map $f_1$ is continuous, there exists $\delta'>0$ such that for all $x$ satisfying $d_{X}(x, x^{0}) < \delta'$, it holds that $d_{X_N}(f_1(x), x^{0}_{N}) < \delta$.
Then for all $x \in X$ such that $d_{X}(x, x^{0}) < \delta'$, it will be true that $d_{X_N}(gf_1(x), gx^{0}_{N}) < \e$ for all $g \in G$.
Thus, for all $g \in G$ and all $x \in X$ satisfying $d_{X}(x, x^{0}) < \delta'$ we have
\begin{multline*}
d_X(gx, gx^{0}) \le \\
\le d_{X_N}(gf_1(x), gx^{0}_{N}) + \bigl| d_X(gx, gx^{0}) - d_{X_N}(gf_1(x), gx^{0}_{N}) \bigr| \le \e + \dis^{G}(f_1) \le 3\e.
\end{multline*}
\end{proof}

Note that in Theorem~\ref{thm:Stability} one cannot replace the convergence with respect to $d^{G}_{GH, p, c}$ by the convergence with respect to $d^{G}_{GH, p}$, as illustrated by the following example.

\begin{examp}
Let $G = \R$. Consider as pointed dynamical systems $(X_n, d_{X_{n}}, x_{n}^{0})$ the subsets
$X_n = \{-1/n\} \cup Q \subset \R$, where $Q = [0, 1] \subset \R$ is a segment with the induced metric.
Let $x_{n}^{0} = -1/n$, and the distinguished point of the limit system $X = Q$ be $x^0 = 0$.

The dynamics on the segment $Q$ (both in $X_{n}$ and in $X$) is defined by the phase flow of the differential equation $\dot{x} = x(1 - x)$ on the line.
The segment $Q$ is invariant since its endpoints $x=0$ and $x=1$ are equilibrium points. The isolated point $x_{n}^{0} = -1/n$ in the system $X_{n}$ is assumed to be fixed.

The systems $(X_n, x^{0}_{n})$ are asymptotically stable. Indeed, for each fixed $n$ the point $x_{n}^{0}$ is isolated in $X_{n}$, therefore it has an open neighborhood containing no other points of the space, and the trajectory does not leave it.

At the same time, $d^{G}_{GH,p}(X_n, X) \to 0$ as $n \to \infty$.
Indeed, consider the maps
$$
f_1 \: X \to X_{n}, \
f_1(x) =
\begin{cases}
-1/n, & \text{ if } x = 0, \\
x, & \text{ if } x \ne 0;
\end{cases}
$$
$$
f_2 \: X_{n} \to X, \
f_2(x) =
\begin{cases}
0, & \text{ if } x = -1/n, \\
x, & \text{ if } x \ne -1/n.
\end{cases}
$$
Then $\dis^{G}(f_1) = \dis^{G}(f_2) = \codis^{G}(f_1, f_2) = 1/n$, which implies that $2d^{G}_{GH,p}(X_{n}, X) \le 1/n$.

However, the limit system $X$ is unstable because the point $x^0$ is no longer isolated, and the condition $\dot{x} > 0$ for $x \in (0, 1)$ implies that for all $x \in (0, 1)$ the trajectories monotonically move away from $x^0$ towards the point $1$.
\end{examp}

\begin{dfn}
We will say that a sequence of quasi-asymptotically stable pointed dynamical systems $(X_{n}, d_{X_{n}}, x^{0}_{n})$ \emph{has a common radius of attraction} $\dl_{0} > 0$, if $\dl_{0}$ is a radius of attraction for all $X_{n}$.
\end{dfn}

\begin{thm}\label{thm:AS}
Let a sequence of quasi-asymptotically stable pointed dynamical systems
$$\{(X_n, d_{X_n}, x^{0}_{n})\}_{n \in \N},$$
having a common radius of attraction $\dl_{0} > 0$, converge with respect to the pointed Gromov--Hausdorff distance to a dynamical system $(X, d_X, x^{0})$\/\rom:
$$d^{G}_{GH, p}\bigl((X_n, x^{0}_{n}), (X, x^{0})\bigr) \to 0.$$
Then the pointed dynamical system $(X, d_{X}, x^{0})$ is quasi-asymptotically stable with the radius of attraction $\dl_0$.
\end{thm}

\begin{proof}
Fix $\e > 0$.
Let $x \in B_{X}(x^{0}, \dl_{0})$.
Let $r = d_{X}(x^{0}, x)$.
Choose $N \in \N$ such that
$$d^{G}_{GH, p}\bigl((X, x^{0}), (X_N, x^{0}_{N})\bigr) \le
\frac{1}{10} \min \{\e, \dl_0 - r\}.$$
Then there exists a map
$f_1 \: X \to X_N, \ f_1(x^{0}) = x^{0}_{N}$ with distortion
$$\dis^{G}(f_1) < \frac{1}{4} \min \{\e, \dl_0 - r\}.$$
Since $\dis^{G}(f_{1}) < (\dl_{0} - r) / 4$, we have
$$f_1(x) \in B_{X_{N}}(x^{0}_{N}, r + \dis^{G}(f_1)) \subset B_{X_{N}}(x^{0}_{N}, r + (\dl_0 - r) / 4) \subset B_{X_{N}}(x^{0}_{N}, \dl_0).$$
Consider a set $B^{*} \in \fB_{G}$ such that for all $g \in B^{*}$ we have $d_{X_{N}}(gf_{1}(x), gx^{0}_{N}) < \e$.
Then for the chosen point $x \in X$ and $g \in B^{*}$ the following holds:
$$d_X(gx, gx^{0}) \le d_{X_N}(gf_1(x), gx^{0}_{N}) +
\bigl|d_X(gx, gx^{0}) - d_{X_N}(gf_1(x), gx^{0}_{N}) \bigr| \le \e + \dis^{G}(f_1) \le 2\e.$$
\end{proof}

Note that in Theorem~\ref{thm:AS} the condition of having a common radius of attraction is very important.
The following example shows that this condition cannot be omitted, even if we additionally assume that, first, all systems $X_n$ are asymptotically stable, second, that $X$ is stable, and third, that the convergence holds with respect to the pointed continuous distance $d^{G}_{GH,p,c}$.

\begin{examp}
Let us show that in Theorem~\ref{thm:AS} the requirement of a common radius of attraction is essential. Consider on the space $\R$ with the distinguished point $x_0 = 0$ the limit dynamical system $X$, defined by the equation:
$$\dot{x} = v(x), \quad v(x) = x^2 \sin\biggl(\frac{1}{x}\biggr), \quad v(0) = 0.$$
Since the function $v(x)$ is continuous and $v(0) = 0$, the point $x_0 = 0$ is an equilibrium point. The stability of $X$ follows from the structure of the neighborhoods of the point $x_0$: in any arbitrarily small neighborhood of $x_0$ to the right and to the left there are zeros $x_{\ell} < 0 < x_{r}$ of the vector field $v$. By the uniqueness of solutions, the trajectory of any point from the interval $(x_{\ell}, x_{r})$ cannot cross its boundaries. Therefore, the interval $(x_{\ell}, x_{r})$ is invariant, which guarantees the Lyapunov stability of the point $x_0 = 0$.

However, the system $X$ is not asymptotically stable. Indeed, in any neighborhood of the point $x_0 = 0$ the function $v(x)$ changes sign infinitely many times. The trajectories cannot overcome these zeros of the field $v(x)$ or ``step over'' them towards $x_0$ , because the zeros themselves are fixed points, and the solution is continuous.

Now let us construct a sequence of approximating systems $X_N$ on $\R$ with the same distinguished point $x_0 = 0$. Let us choose a sequence of points $A_N \to 0$ such that the sine evaluated at them takes the value $-1$, for example, $A_N = \frac{1}{3\pi/2 + 2\pi N}$. Then $v(A_N) = -A_N^2 < 0$. We define the systems $X_N$ by the equations $\dot{x} = v_N(x)$, where the vector field $v_N$ is obtained by ``cutting off'' the original field on the segment $[-A_N, A_N]$ and replacing it with a linear function:
$$
v_N(x) =
\begin{cases}
    v(x), & \text{for } |x| \ge A_N, \\
    -A_N x, & \text{for } |x| \le A_N.
\end{cases}
$$
For each fixed $N$ the system $X_N$ in a neighborhood of zero coincides with the linear stable system $\dot{x} = -A_N x$, which trivially implies its asymptotic stability with the radius of attraction $\delta_N = A_N$.

Choosing the identity maps as $f_1$ and $f_2$, we obtain
$$d^{G}_{GH,p,c}(X_N, X) \to 0 \text{ as } N \to \infty.$$ However, the systems $X_N$ do not have a common radius of attraction, since $\delta_N = A_N \to 0$, therefore the conditions of Theorem~\ref{thm:AS} are not satisfied here.
\end{examp}

From Theorems~\ref{thm:Stability} and~\ref{thm:AS}, we obtain the following corollary.

\begin{cor}\label{cor:AS}
Let a sequence of asymptotically stable pointed dynamical systems $\{(X_n, d_{X_n}, x^{0}_{n})\}_{n \in \N}$ having a common radius of attraction $\dl_{0} > 0$, converge with respect to the pointed continuous Gromov-Hausdorff distance to a dynamical system $(X, d_X, x^{0})$\/\rom:
$$d_{GH, p, c}\bigl((X_n, x^{0}_{n}), (X, x^{0})\bigr) \to 0.$$
Then the pointed dynamical system $(X, d_{X}, x^{0})$ is asymptotically stable; moreover, $\dl_{0}$ is its radius of attraction.
\end{cor}

\subsection{Compact quasi-asymptotic stability}

\begin{dfn}
Let $(X, d_{X}, \phi_{X})$ be a dynamical system.
We say that a set $M \subset X$ \emph{is uniformly attracted} to a point $x_0 \in X$ (along a base $\fB_{G}$), if the distance from the images of the points of this set to $x_0$ tends uniformly to zero:
$$\lim_{\fB_{G}} \sup_{x \in M} d_{X}(gx, gx_0) = 0.$$
In other words, for any $\e > 0$ there exists an element $B \in \fB_{G}$ such that for all $g \in B$ and for all $x \in M$ we have $d_{X}(gx, gx_0) < \e$.
\end{dfn}

\begin{dfn}\label{dfn:LocUniComp}
Let $X$ be a dynamical system. We will say that a point $x_{0} \in X$ \emph{locally uniformly attracts compact sets} (or that the class of compact sets \emph{is locally uniformly attracted} to the point $x_0$), if there exists an open neighborhood $V \subset X$ of the point $x_0$ such that any compact set $K \subset V$ contained in it is uniformly attracted to the point $x_0$.
\end{dfn}

\begin{dfn}
A pointed dynamical system $(X, x_0)$ \emph{is called compactly quasi-asymptotically stable} if the point $x_0$ locally uniformly attracts compact sets.
\end{dfn}

Just as before, we will say that a number $\dl > 0$ is a radius of attraction for a compactly quasi-asymptotically stable pointed dynamical system $X$, if in Definition~\ref{dfn:LocUniComp} as $V$ one can take the open ball $B_{X}(x_0, \dl)$.

\begin{thm}\label{thm:compactAS}
Let a sequence of compactly quasi-asymptotically stable pointed dynamical systems $\{(X_n, d_{X_n}, x^{0}_{n})\}_{n \in \N}$ having a common radius of attraction $\dl_{0} > 0$, converge with respect to the pointed continuous Gromov--Hausdorff distance to a dynamical system $(X, d_X, x^{0})$\/\rom:
$$d^{G}_{GH, p, c}\bigl((X_n, x^{0}_{n}), (X, x^{0})\bigr) \to 0.$$
Then the pointed dynamical system $(X, d_{X}, x^{0})$ is compactly quasi-asymptotically stable with the radius of attraction $\dl_0$.
\end{thm}

\begin{proof}
Fix $\e > 0$.
Let $K \subset B_{X}(x^{0}, \dl_{0})$.
Let $r = \sup_{x \in K}d_{X}(x, x^{0}) < \dl_{0}$.
Choose $N \in \N$ such that
$$d^{G}_{GH, p, c}\bigl((X, x^{0}), (X_N, x^{0}_{N})\bigr) \le
\frac{1}{10} \min \{\e, \dl_0 - r\}.$$
Then there exists a continuous map
$f_1 \: X \to X_N, \ f_1(x^{0}) = x^{0}_{N}$ with distortion
$$\dis^{G}(f_1) < \frac{1}{4} \min \{\e, \dl_0 - r\}.$$
Note that $K' = f_1(K)$ is compact as the image of a compact set under a continuous map.
Since $\dis^{G}(f_{1}) < (\dl_{0} - r)/ 4$, we have
$$K' = f_1(K) \subset B_{X_{N}}(x^{0}_{N}, r + \dis^{G}(f_1)) \subset B_{X_{N}}(x^{0}_{N}, r + (\dl_0 - r) / 4) \subset B_{X_{N}}(x^{0}_{N}, \dl_0).$$
Consider $B^{*} \in \fB$ such that for all $g \in B^{*}$ and all $x \in K'$ we have $d_{X_{N}}(gf_{1}(x), gx^{0}_{N}) < \e$.
Then for all points $x \in K$ and $g \in B^{*}$ the following holds:
$$d_X(gx, gx^{0}) \le d_{X_N}(gf_1(x), gx^{0}_{N}) +
\bigl|d_X(gx, gx^{0}) - d_{X_N}(gf_1(x), gx^{0}_{N}) \bigr| \le \e + \dis^{G}(f_1) \le 2\e.$$
\end{proof}

\begin{dfn}
A pointed dynamical system is called \emph{compactly asymptotically stable} if it is stable and compactly quasi-asymptotically stable.
\end{dfn}

From Theorems~\ref{thm:Stability} and~\ref{thm:compactAS} we obtain the following corollary.

\begin{cor}\label{cor:compactAS}
Let a sequence of compactly asymptotically stable pointed dynamical systems $\{(X_n, d_{X_n}, x^{0}_{n})\}_{n \in \N}$ having a common radius of attraction $\dl_{0} > 0$, converge with respect to the pointed continuous Gromov--Hausdorff distance to a dynamical system $(X, d_X, x^{0})$\/\rom:
$$d^{G}_{GH, p, c}\bigl((X_n, x^{0}_{n}), (X, x^{0})\bigr) \to 0.$$
Then the pointed dynamical system $(X, d_{X}, x^{0})$ is compactly asymptotically stable; moreover, $\dl_{0}$ is its radius of attraction.
\end{cor}

\subsection{Nowhere density of the space of stable dynamical systems}

\begin{dfn}
Let $X, Y$ be dynamical systems. Let $1 \le p < \infty$.
We call the $l^p$-\emph{product of dynamical systems} $X \times_{l^{p}} Y$ the dynamical system $X \times Y$ with the map
$$\phi_{X \times Y} \: G \times X \times Y \to X \times Y, \
g(x, y) = (gx, gy)$$
and the generalized pseudometric
$$d_{X \times Y}\bigl((x_1, y_1), (x_2, y_2)\bigr) =
\begin{cases}
\bigl(d_{X}(x_1, x_2)^{p} + d_{Y}(y_1, y_2)^{p}\bigr)^{1/p}, & \text{if } d_{X}(x_1, x_2), d_{Y}(y_1, y_2) < \infty, \\
\infty, & \text{otherwise}.
\end{cases}$$
\end{dfn}

\begin{lem}\label{lem:bq}
Let $(X, x_0), (Y, y_0)$ be pointed dynamical systems.
Then $$2d^{G}_{GH, p, c}(X, X \times_{l^{1}} Y) \le \diam(Y).$$
\end{lem}

\begin{proof}
Consider the continuous maps $f_1\: X \to X \times Y, \ x \mapsto (x, y_0)$ --- a natural isometric embedding, $\dis^{G}(f_1) = 0$ and
$f_2 \: X \times Y \to X, \ (x, y) \mapsto x$ --- the projection onto $X$.
Let us compute the distortion of $f_2$:
\begin{multline*}
\dis^{G}(f_2) = \sup_{g \in G} \sup_{\substack{x_1, x_2 \in X \\ y_1, y_2 \in Y}} \Bigl| d_{X \times Y}\bigl(g(x_1, y_1), g(x_2, y_2)\bigr) - d_{X}(gx_1, gx_2)\Bigr| = \\
= \sup_{g \in G} \sup_{\substack{x_1, x_2 \in X \\ y_1, y_2 \in Y}} \Bigl| d_{X}(gx_1, gx_2) + d_{Y}(gy_1, gy_2) - d_{X}(gx_1, gx_2)\Bigr|
= \sup_{g \in G} \sup_{y_1, y_2 \in Y} d_{Y}(gy_1, gy_2) = \diam(Y).
\end{multline*}
Let us compute the codistortion:
\begin{multline*}
\codis^{G}(f_1, f_2) = \sup_{g \in G} \ \sup_{\substack{x_1 \in X \\ w = (x_2, y_2) \in X \times Y}} \
\Bigl| d_{X}\bigl(gx_1, gf_2(w)\bigr) - d_{X \times Y}\bigl(gw, gf_1(x_1)\bigr) \Bigr| = \\
= \sup_{g \in G} \sup_{x_1 \in X} \sup_{\substack{x_2 \in X \\ y_2 \in Y}}
\bigl|d_{X}(gx_{1}, gx_{2}) - d_{X}(gx_2, gx_1) - d_{Y}(gy_2, gy_0)\bigr| = \diam(Y).
\end{multline*}
Therefore,
$2d^{G}_{GH, p, c}(X, X \times_{l^{1}} Y) \le \max \bigl\{ \dis^{G}(f_1), \ \dis^{G}(f_2), \ \codis^{G}(f_1, f_2) \bigr\} = \diam(Y)$.
\end{proof}

\begin{thm}\label{thm:dns}
Let $(G, e)$ be an arbitrary pointed metric space containing at least one element other than $e$. Then the subset of all stable pointed dynamical systems is closed and nowhere dense in the space $\cDM_{p, c}$.
\end{thm}

\begin{proof}
Denote by $\cS$ the subset of all stable systems. From the previously proved Theorem~\ref{thm:Stability} it follows that the set $\cS$ is closed in the topology of the continuous Gromov--Hausdorff distance. It is a well-known fact that a closed set is nowhere dense if and only if every neighborhood of any of its points contains a point from its complement. Thus, it is sufficient to show that in any neighborhood of an arbitrary stable system there is an unstable one.

Fix a stable system $(X, d_X, x^0)$ and an arbitrary $\e > 0$.
Consider the segment $I_{\e} = ([0, \e], 0)$ as a pointed dynamical system with the defined map
$$\phi_{I_{\e}} \: G \times I_{\e} \to I_{\e}, \
(g, x) \mapsto gx =
\begin{cases}
x, & \text{if } g = e, \\
\e, & \text{if } g \ne e \text{ and } x \ne 0, \\
0, & \text{if } x = 0.
\end{cases}
$$
Now let $X_{\e} = X \times_{l^1} I_{\e}$; this system is unstable because in a neighborhood of any point $(x_0, 0)$ there are points distinct from it in the segment $\{x\} \times I_{\e}$, which under the action of any element $g \ne e$ are mapped to the point $(gx, \e)$.

By Lemma~\ref{lem:bq} we get
$d^{G}_{GH, p, c}(X, X_{\e}) \le \e / 2$.
Thus, we have shown that arbitrarily close to the system $X$ there is an unstable system $X_{\e}$. Therefore, the set $\cS$ is nowhere dense in the space $\cDM_{p, c}$.
\end{proof}

\begin{thm}
Let $(G, e)$ be an arbitrary nontrivial pointed metric space. Then
\begin{enumerate}
\item
the subspace of all asymptotically stable compact dynamical systems is a nowhere dense $F_{\sigma}$-subset of $\cDM_{p, c}$;
\item
the subspace of all compactly asymptotically stable compact dynamical systems is a nowhere dense $F_{\sigma}$-subset of $\cDM_{p, c}$.
\end{enumerate}
\end{thm}

\begin{proof}
We will prove both statements simultaneously.
Denote by $S$ the subspace of interest.
From Theorem~\ref{thm:dns} it follows that $S$ is nowhere dense in $\cDM_{p, c}$ as a subset of a nowhere dense set.
For any $n \in \N$ denote by $S_n \subset S$ the subset of all (compactly) asymptotically stable dynamical systems with radius of attraction $1/n$.
By Corollary~\ref{cor:AS} or~\ref{cor:compactAS} we obtain that
$S_{n}$ is closed in $\cDM_{p, c}$.
Consequently, their union
$S = \cup_{n = 1}^{\infty} S_{n}$ is an $F_{\sigma}$-subset of $\cDM_{p,c}$.
\end{proof}

\subsection{Quasi-asymptotic stability with respect to a class of subsets}

In this section, we introduce the concept of quasi-asymptotic stability with respect to an arbitrary class of subsets and generalize Theorem~\ref{thm:compactAS} to this case.

\begin{notation}
For any set $X$ denote by $\cP_{0}(X)$ the set of all non-empty subsets of $X$.
\end{notation}

\begin{notation}
For any generalized pseudometric space $X$ denote by $\cP_{0}^{\text{comp}}(X)$ the class of all non-empty compact subsets of $X$.
\end{notation}

\begin{dfn}\label{dfn:LocUni}
Let $X$ be a dynamical system, and $\fU \subset \cP_{0}(X)$ be an arbitrary subclass of non-empty subsets of $X$.
We will say that the class $\fU$ \emph{is locally uniformly attracted} to the point $x_0$ (or that the point $x_0$ \emph{locally uniformly attracts the class} $\fU$), if there exists an open neighborhood $V \subset X$ of the point $x_{0}$ such that any set $U \in \fU$ satisfying the condition $U \subset V$, is uniformly attracted to the point $x_0$.
\end{dfn}

\begin{dfn}
We will say that a pointed dynamical system $(X, x_0)$ \emph{is quasi-asymptotically stable} with respect to the class $\fU \subset \cP_{0}(X)$, if this class $\fU$ is locally uniformly attracted to the point $x_{0}$.
\end{dfn}

\begin{rk}
It follows from the definition that a pointed dynamical system $(X, x_{0})$ is compactly quasi-asymptotically stable if and only if it is quasi-asymptotically stable with respect to the class $\cP^{\text{comp}}_{0}(X)$ of all non-empty compact subsets of $X$.
\end{rk}

Let $X$ be a pointed dynamical system.
Consider $\bU_{0}(X)$, the family of all classes of subsets $\fU \subset \cP_{0}(X)$ that are locally uniformly attracted to the point $x_0$ in the sense of the previous definition.

\begin{ass}
The family $\bU_0(X)$ possesses the following properties:
\begin{enumerate}
\item \label{rk:U0:1} $\bU_0(X) \ne \emptyset$;
\item \label{rk:U0:2} if $\fU_1 \in \bU_0(X)$ and $\fU_2 \subset \fU_1$, then $\fU_2 \in \bU_{0}(X)$;
\item \label{rk:U0:3} if $\fU_1 \in \bU_0(X)$ and $\fU_2 \in \bU_0(X)$, then $\fU_1 \cup \fU_2 \in \bU_{0}(X)$;
\item \label{rk:U0:4} if $\fU \in \cP_{0}(X), \ \fU' \in \bU_0(X)$ and for any $U \in \fU$ there exists $U' \in \fU'$ such that $U \subset U'$, then $\fU' \in \bU_0(X)$;
\item \label{rk:U0:5} if $\fU_1, \fU_2 \in \bU_0(X)$, then $\{U_1 \cup U_2 \mid U_1 \in \fU_1, \ U_2 \in \fU_2 \} \in \bU_{0}(X)$.
\end{enumerate}
\end{ass}

\begin{proof}
Properties~\ref{rk:U0:1},~\ref{rk:U0:2}, and~\ref{rk:U0:4} follow from the definition. It remains to check properties~\ref{rk:U0:3} and~\ref{rk:U0:5}.

Let us check property~\ref{rk:U0:3}. If $V_1, V_2$ are neighborhoods from Definition~\ref{dfn:LocUni} for $\fU_1$ and $\fU_2$ respectively, then for $\fU_1 \cup \fU_2$ it is sufficient to take the neighborhood $V_1 \cap V_2$.

Let us check property~\ref{rk:U0:5}. Consider the same neighborhood $V$, fix $\e > 0$.
Then for any $U_1 \in \fU_1, \ U_1 \subset V$ there exists $B_{1}^{*} \in \fB_{G}$ such that for all $x \in U_1$ and $g \in B_{1}^{*}$ we have $d_{X}(gx, gx_0) < \e$. For any $U_2 \in \fU_2, \ U_2 \subset V$ there exists $B_{2}^{*} \in \fB_{G}$ such that for all $x \in U_2$ and $g \in B_{2}^{*}$ we have $d_{X}(gx, gx_0) < \e$. Therefore, for any $B^{*} \in \fB_{G}$ satisfying $B^{*} \subset B_1^{*} \cap B_2^{*}$, and for all $x \in U_1 \cup U_2$ and $g \in B^{*}$ it will hold that $d_{X}(gx, gx_0) < \e$.
\end{proof}

Now let us generalize Theorem~\ref{thm:compactAS} to the case of an arbitrary class $\fU \in \bU_0$.
Let each pointed dynamical system $X$ be associated with its own class of subsets $\fU(X) \in \bU_{0}(X)$.

\begin{thm}
Let a sequence of pointed dynamical systems $$\{(X_n, d_{X_n}, x_{n}^{0})\}_{n \in \N}$$ converge with respect to the categorical Gromov--Hausdorff distance $d^{G}_{GH, \cK}$ to a pointed dynamical system $(X, d_{X}, x^{0})$\/\rom:
$$d^{G}_{GH,\cK}(X_{n}, X) \to 0.$$
Suppose that for all $n \in \N$ the system $X_n$ is quasi-asymptotically stable with respect to the class $\fU(X_n)$, and the number $\dl_{0} > 0$ is a common radius of attraction for them.
Also suppose that for any $f \in \Mor_{\cK}(X, X_n)$ and any $U \in \fU(X)$ we have $f(U) \in \fU(X_n)$.
Then the pointed dynamical system $X$ is quasi-asymptotically stable with respect to the class $\fU(X)$; moreover, any number $0 < r < \dl_0$ is its radius of attraction.
\end{thm}

\begin{proof}
The proof conceptually and technically completely follows the scheme of the proof of Theorem~\ref{thm:compactAS}.
\end{proof}

\section{Hausdorff map}

\subsection{The Hausdorff map is 1-Lipschitz}

\begin{dfn}
Let $(X, d_{X})$ be a generalized pseudometric space.
For non-empty $A, B \subset X$ we define the \emph{Hausdorff distance} between $A$ and $B$:
$$d^{H}_{d_{X}}(A, B) = \max \biggl\{ \sup_{a \in A} \inf_{b \in B} d_{X}(a,b), \
\sup_{b \in B} \inf_{a \in A} d_{X}(a, b) \biggr\}.$$
\end{dfn}

\begin{notation}
For an arbitrary metric space $X$ denote by $\cH(X)$ the family of all non-empty closed bounded subsets of $X$.
\end{notation}

\begin{prop}[\cite{BurBurIva}] \label{prop:H}
For every metric space $X$ the Hausdorff distance is a metric on $\cH(X)$.
Moreover, the spaces $X$ and $\cH(X)$ simultaneously possess any of the following properties\/\rom: completeness, boundedness, total boundedness, compactness.
\end{prop}

\begin{thm}[\cite{Mich}]
For any non-empty compact metric spaces $X$ and $Y$ it holds that
$$d_{GH}\bigl(\cH(X), \cH(Y) \bigr) \le d_{GH}(X, Y).$$
\end{thm}

Let us extend this theorem to the case of $F$-spaces.
For an arbitrary $F$-space
$$\bigl(X, \{d_{X, i}\}_{i \in I}, \{x_{j}\}_{j \in J}\bigr)$$ on the space of all its non-empty subsets $\cP_{0}(X)$ there naturally arises an $F$-space structure with the family of generalized pseudometrics
$$d_{\cP_{0}(X), i}(A, B) = d^{H}_{d_{X, i}}(A, B)$$
and distinguished points $\bigl\{\{x_{j}\}\bigr\}_{j \in J}$.

We define the \emph{Hausdorff functor} as the functor
$$\cP_{0} \: \cF \to \cF, \ X \mapsto \cP_{0}(X), \quad \cP_{0} \: \Mor_{\cF}(X, Y) \to \Mor_{\cF}\bigl(\cP_{0}(X), \cP_{0}(Y)\bigr),$$
which associates to a morphism $f \in \Mor_{\cF}(X, Y)$ the morphism
$$F \in \Mor_{\cF}\bigl(\cP_{0}(X), \cP_{0}(Y)\bigr),$$
defined on any set $A \in \cP_{0}(X)$ as the image under $f$, that is, $F(A) := f(A)$.

\begin{lem}\label{lem:matan}
Let $a, b, c, d \in \R$. Then
$$\bigl| \max\{a, b\} - \max\{c, d\}\bigr| \le \max\bigl\{ |a - c|, |b - d|\bigr\}.$$
\end{lem}

\begin{proof}
Note that $a = c + (a - c) \le c + |a - c| \le \max\{c, d\} + \max\bigl\{|a - c|, |b - d|\bigr\}$.
Similarly, $b = d + (b - d) \le d + |b - d| \le \max\{c, d\} + \max\bigl\{|a - c|, |b - d|\bigr\}$.
Therefore,
$$ \max\{a, b\} \le \max\{c, d\} + \max\bigl\{|a - c|, |b - d|\bigr\}. $$
Swapping the roles of the pairs $(a, b)$ and $(c, d)$, by symmetry we obtain a similar inequality for $\max\{c, d\}$. Combining them, we arrive at the desired statement.
\end{proof}

\begin{lem}\label{lem:supinf}
Let $A$ and $B$ be arbitrary non-empty sets, and $F, G \: A \times B \to \R$ be bounded functions. Then
$$ \Bigl| \sup_{a \in A} \inf_{b \in B} F(a,b) - \sup_{a \in A} \inf_{b \in B} G(a,b) \Bigr| \le \sup_{a \in A, b \in B} \bigl| F(a,b) - G(a,b) \bigr|. $$
\end{lem}

\begin{proof}
For any fixed $a \in A$ and $b \in B$ the following obvious inequality holds
$$F(a,b) \le G(a,b) + |F(a,b) - G(a,b)|.$$
Taking the infimum over $b$, we obtain:
$$ \inf_{b \in B} F(a,b) \le \inf_{b \in B} G(a,b) + \sup_{b \in B} |F(a,b) - G(a,b)|. $$
Now taking the supremum over $a$:
$$\sup_{a \in A} \inf_{b \in B} F(a,b)
\le \sup_{a \in A} \inf_{b \in B} G(a,b) + \sup_{a \in A} \sup_{b \in B} |F(a,b) - G(a,b)|.$$
It follows that
$$ \sup_{a \in A} \inf_{b \in B} F(a,b) - \sup_{a \in A} \inf_{b \in B} G(a,b) \le \sup_{a \in A, b \in B} |F(a,b) - G(a,b)|. $$
Swapping the functions $F$ and $G$ , by symmetry we obtain the required estimate.
\end{proof}

\begin{thm}\label{thm:Hausdorff}
Suppose that the subcategories $\cK, \cK' \subset \cF$ are chosen such that for all $X, Y \in \Ob(\cK)$ we have
$$\cP_{0}(\Mor_{\cK}(X, Y)) \subset \Mor_{\cK'}( \cP_{0}(X), \cP_{0}(Y)),$$
that is, the restriction of the Hausdorff functor $\cP_{0}$ to $\cK$ is a well-defined functor $\cP_{0}|_{\cK} \: \cK \to \cK'$.
Then for any $X, Y \in \Ob(\cK)$ it holds that
$$d^{F}_{GH, \cK'}\bigl(\cP_{0}(X), \cP_{0}(Y) \bigr) \le d^{F}_{GH, \cK}(X, Y).$$
\end{thm}

\begin{proof}
We can assume that $d_{GH, \cK}(X, Y) < \infty$, since otherwise there is nothing to prove.
Fix some $\e > 0$ and consider the maps
$f_1 \in \Mor_{\cK}(X, Y), \ f_2 \in \Mor_{\cK}(Y, X)$, for which
$$\max\bigl\{ \dis(f_1), \dis(f_2), \codis(f_1, f_2)\bigr\} \le 2d_{GH, \cK}(X, Y) + \e.$$

Let $F_1 = \cP_{0}(f_1), \ F_2 = \cP_{0}(f_2)$.
Let us estimate the distortion of the map $F_1$. Thus,
\begin{multline*}
\dis(F_1) =
\sup_{X_1, X_2 \in \cP_{0}(X)} \ \sup_{i \in I} \
\Bigl|
d_{\cP_{0}(X), i}(X_1, X_2) - d_{\cP_{0}(Y), i}\bigl(f_1(X_1), f_1(X_2)\bigr)
\Bigr| = \\
\
= \sup_{X_1, X_2 \in \cP_{0}(X)} \ \sup_{i \in I} \
\Biggl|
\max \biggl\{
\sup_{x_1 \in X_1} \inf_{x_2 \in X_2} d_{X, i}(x_1, x_2), \
\sup_{x_2 \in X_2} \inf_{x_1 \in X_1} d_{X, i}(x_1, x_2)
\biggr\} \ - \\
\
- \max \biggl\{
\sup_{x_1 \in X_1} \inf_{x_2 \in X_2} d_{Y, i}\bigl(f_1(x_1), f_1(x_2)\bigr), \
\sup_{x_2 \in X_2} \inf_{x_1 \in X_1} d_{Y, i}\bigl(f_1(x_1), f_1(x_2)\bigr)
\biggr\}
\Biggr| \le \\
\
\le \sup_{X_1, X_2 \in \cP_{0}(X)} \ \sup_{i \in I} \
\max \Biggl\{
\biggl| \sup_{x_1 \in X_1} \inf_{x_2 \in X_2} d_{X, i}(x_1, x_2) - \\ - \sup_{x_1 \in X_1} \inf_{x_2 \in X_2} d_{Y, i}\bigl(f_1(x_1), f_1(x_2)\bigr) \biggr|, \\
\
\biggl| \sup_{x_2 \in X_2} \inf_{x_1 \in X_1} d_{X, i}(x_1, x_2) - \sup_{x_2 \in X_2} \inf_{x_1 \in X_1} d_{Y, i}\bigl(f_1(x_1), f_1(x_2)\bigr) \biggr|
\Biggr\} \le
\end{multline*}
\begin{multline*}
\le \sup_{X_1, X_2 \in \cP_{0}(X)} \ \sup_{i \in I} \ \max \Biggl\{
\sup_{x_1 \in X_1, x_2 \in X_2} \Bigl| d_{X, i}(x_1, x_2) - d_{Y, i}\bigl(f_1(x_1), f_1(x_2)\bigr) \Bigr|, \\
\
\sup_{x_1 \in X_1, x_2 \in X_2} \Bigl| d_{X, i}(x_1, x_2) - d_{Y, i}\bigl(f_1(x_1), f_1(x_2)\bigr) \Bigr|
\Biggr\} = \\
\
= \sup_{X_1, X_2 \in \cP_{0}(X)} \ \sup_{i \in I} \ \sup_{x_1 \in X_1, x_2 \in X_2} \Bigl| d_{X, i}(x_1, x_2) - d_{Y, i}\bigl(f_1(x_1), f_1(x_2)\bigr) \Bigr| \le \\
\
\le \sup_{x_1, x_2 \in X} \ \sup_{i \in I} \ \Bigl| d_{X, i}(x_1, x_2) - d_{Y, i}\bigl(f_1(x_1), f_1(x_2)\bigr) \Bigr| = \dis(f_1).
\end{multline*}
Similarly, we obtain that $\dis(F_2) \le \dis(f_2)$.

In the same way, we estimate the codistortion:
\begin{multline*}
\codis(F_1, F_2) = \sup_{\substack{X' \in \cP_{0}(X) \\ Y' \in \cP_{0}(Y)}} \sup_{i \in I} \
\biggl|
d_{\cP_{0}(X), i}\bigl(X', F_2(Y') \bigr)
- d_{\cP_{0}(Y), i}\bigl(Y', F_1(X')\bigr)
\biggr| \le \\
\
\le \sup_{\substack{X' \in \cP_{0}(X) \\ Y' \in \cP_{0}(Y)}} \sup_{i \in I} \sup_{\substack{x \in X' \\ y \in Y'}} \
\Bigl| d_{X, i}\bigl(x, f_2(y)\bigr) - d_{Y, i}\bigl(y, f_1(x)\bigr) \Bigr| \le \\
\
\le \sup_{\substack{x \in X \\ y \in Y}} \sup_{i \in I} \
\Bigl| d_{X, i}\bigl(x, f_2(y)\bigr) - d_{Y, i}\bigl(y, f_1(x)\bigr) \Bigr| = \codis(f_1, f_2).
\end{multline*}

Thus,
\begin{multline*}
2d_{GH, \cK'}(\cP_{0}(X), \cP_{0}(Y)) \le \max \bigl\{ \dis(F_1), \dis(F_2), \codis(F_1, F_2)\bigr\} \le \\
\le \max \bigl\{ \dis(f_1), \dis(f_2), \codis(f_1, f_2)\bigr\} \le
2d_{GH, \cK}(X, Y) + \e.
\end{multline*}
Letting $\e \to 0$, we obtain the desired result.
\end{proof}

\subsection{Hausdorff map for continuous dynamical systems}

Note that in the definition of an $F$-space and a dynamical system, we considered generalized pseudometrics.

\begin{dfn}
We define a \emph{metric dynamical system} as a dynamical system $(X, d_{X})$ such that $d_{X}$ is a metric.
\end{dfn}

\begin{notation}
For an arbitrary topological space $X$ and a subset $A \subset X$ we will denote by $\overline{A}$ its closure.
\end{notation}

\begin{rk}
Note that for any non-empty subsets $A, B \subset X$ of a metric space $(X, d_{X})$ the Hausdorff distance between them is equal to the Hausdorff distance between their closures:
$$d^{H}_{d_{X}}(A, B) = d^{H}_{d_{X}}(\overline{A}, \overline{B}).$$
\end{rk}

\begin{dfn}
We will say that a dynamical system $X$
\emph{is uniformly continuous} if for all $g \in G$ the map
$\phi_{X}(g, \bullet) \: X \to X$ is uniformly continuous.
\end{dfn}
\begin{dfn}
We will say that a dynamical system $X$ is endowed with a
\emph{bounded action} if for all $g \in G$ the map
$\phi_{X}(g, \bullet) \: X \to X$ is bounded, that is, the image of any bounded set under this map is bounded.
\end{dfn}

Denote by $\cDGH_{0} \subset \cDGH$ the subclass of all uniformly continuous metric dynamical systems with a bounded action.

\begin{rk}
Note that the uniform continuity of a map $f \: X \to Y$ (where $X, Y$ are metric spaces) does not, in general, imply its boundedness.
Indeed, consider $X = \R$ with the discrete metric (i.e., for $x \ne y$ by definition $d_{X}(x, y) = 1$), and let $Y$ be $\R$ with the standard metric.
Then the identity map $\id \: X \to Y$ is uniformly continuous, but not bounded.
\end{rk}

For any dynamical system $X \in \cDGH_{0}$ , by analogy with metric spaces, we can also introduce the Hausdorff space $\cH(X)$, consisting of all non-empty closed bounded subsets endowed with the Hausdorff metric.
The dynamical structure is introduced by analogy with $F$-spaces, but now instead of the image, we consider its closure. Namely, given the map $\phi_{X} \: G \times X \to X$, we define the map $\phi_{\cH(X)} \: G \times \cH(X) \to \cH(X), \ A \mapsto \overline{\phi_{X}(g, A)}$. Thus, we have obtained a dynamical system $\cH(X)$.

The following statement immediately follows from the definition of $\cDGH_{0}$, the map $\cH$ and Proposition~\ref{prop:H}.

\begin{ass}\label{ass:CompactContinuous}
If $X \in \cDGH_{0}$, then $\cH(X) \in \cDGH_{0}$.
\end{ass}

\begin{thm}\label{prop:Hausdorff}
Suppose we are given categories $\cK, \cK' \subset \cF$ and $\cDGH_{0} \subset \Ob(\cK) \cap \Ob(\cK')$.
Suppose that for any $X, Y \in \cDGH_{0}$ and any $f \in \Mor_{\cK}(X, Y)$, the map $F \: \cH(X) \to \cH(Y)$ defined by $F(A) := \overline{f(A)}$ is an element of $\Mor_{\cK'}(\cH(X), \cH(Y))$.
Then for any $X, Y \in \cDGH_{0}$ we have
$$d^{G}_{GH, \cK'}\bigl(\cH(X), \cH(Y)\bigr) \le d^{G}_{GH, \cK}(X, Y).$$
\end{thm}

\begin{proof}
The proof will be similar to the proof of Theorem~\ref{thm:Hausdorff}.
We can assume that $d_{GH, \cK}(X, Y) < \infty$.
Fix $\e > 0$ and consider $f_1 \in \Mor_{\cK}(X, Y), \ f_2 \in \Mor_{\cK}(Y, X)$, for which
$$\max\bigl\{ \dis(f_1), \dis(f_2), \codis(f_1, f_2)\bigr\} \le 2d_{GH, \cK}(X, Y) + \e.$$
Let $F_1 = \cH(f_1), \ F_2 = \cH(f_2)$.

Since the dynamical system $X$ is continuous, the action
$\phi_{X}(g, \bullet) \: x \mapsto gx$ maps the closure of a set into a subset of the closure of its image. This implies the identity (for any $A \in \cH(X)$):
$$\phi_{\cH(Y)}\bigl(g, F_1(A)\bigr)
= \overline{\phi_{Y}\bigl(g, F_1(A)\bigr)}
= \overline{\phi_{Y}\bigl(g, \overline{f_1(A)} \bigr)}
= \overline{\phi_{Y}\bigl(g, f_1(A) \bigr)}.$$
A similar identity holds for $\phi_{\cH(X)}$.
Moreover, the Hausdorff distance between sets coincides with the distance between their closures.
Taking this into account, we can now repeat verbatim the entire chain of inequalities from the proof of Theorem~\ref{thm:Hausdorff} to obtain that
$$\dis(F_1) \le \dis(f_1), \ \dis(F_2) \le \dis(f_2), \ \codis(F_1, F_2) \le \codis(f_1, f_2).$$
Therefore,
\begin{multline*}
2d_{GH, \cK'}(\cH(X), \cH(Y)) \le \max \bigl\{ \dis(F_1), \dis(F_2), \codis(F_1, F_2)\bigr\} \le \\
\le \max \bigl\{ \dis(f_1), \dis(f_2), \codis(f_1, f_2)\bigr\} \le
2d_{GH, \cK}(X, Y) + \e.
\end{multline*}
Letting $\e \to 0$, we obtain the desired result.
\end{proof}

\begin{cor}
Let $X, Y$ be arbitrary \rom(not necessarily compact\rom) metric spaces, and let $\cH(X)$ and $\cH(Y)$ be the spaces of their nonempty closed bounded subsets equipped with the Hausdorff metric.
Then
$$d_{GH}\bigl(\cH(X), \cH(Y)\bigr) \le d_{GH}(X, Y),$$
where $d_{GH}$ is the classical Gromov--Hausdorff distance \rom(the infimum is taken over all maps, not necessarily continuous ones\rom).
\end{cor}

\begin{proof}
Any metric space can be considered as a dynamical system with the trivial action $\phi(e, x) = x$.
Such an action is bounded and uniformly continuous, hence $X, Y \in \cDGH_{0}$.
Since in the classical Gromov--Hausdorff distance the classes of morphisms consist of all maps ($\cK = \cK' = \cF$), the condition of Theorem~\ref{prop:Hausdorff} is satisfied automatically: for any map $f \: X \to Y$ the map $F \: A \mapsto \overline{f(A)}$ is a well-defined map from $\cH(X)$ to $\cH(Y)$.
Applying Theorem~\ref{prop:Hausdorff}, we obtain the desired inequality.
\end{proof}

\begin{notation}
For dynamical systems $X, Y \in \cDGH_{0}$ consider the Gromov--Hausdorff distance $d^{G}_{GH_{0}}(X, Y)$, in which as $\cK$ we consider $\cDGH_{0}$, and the morphisms $\Mor_{\cK}(X, Y)$ are all uniformly continuous bounded maps.
\end{notation}

From the proved theorem, we obtain the following corollary.

\begin{cor}\label{cor:Hausdorff_1}
For any $X, Y \in \cDGH_{0}$ it holds that
$$d^{G}_{GH_{0}}(\cH(X), \cH(Y)) \le d^{G}_{GH_{0}}(X, Y).$$
\end{cor}

\begin{proof}
Indeed, it is only necessary to check that the map $\cH$ is well-defined.
This condition is satisfied, since a uniformly continuous bounded map $f \: X \to Y$ defines a uniformly continuous bounded map $F \: \cH(X) \to \cH(Y), \ F(A) = f(A)$.
\end{proof}

\begin{rk}\label{rk:Hausdorff_2}
Note that we have the inclusion $\cDM \subset \cDGH_{0}$.
Also note that if $X, Y \in \cDM$, then any continuous map $f \: X \to Y$ is uniformly continuous and bounded; therefore, in this case the just introduced distance $d^{G}_{GH_{0}}$ coincides with the continuous distance:
$$d^{G}_{GH_{0}}(X, Y) = d^{G}_{GH, c}(X, Y).$$
\end{rk}

From Corollary~\ref{cor:Hausdorff_1} and Remark~\ref{rk:Hausdorff_2}, we obtain the following corollary.

\begin{cor}
Let $X, Y \in \cM$ be non-empty compact metric spaces.
Then
$$d_{GH, c}\bigl(\cH(X), \cH(Y)\bigr) \le d_{GH, c}(X, Y).$$
\end{cor}

Similarly, we obtain the following corollary.

\begin{cor}
\begin{enumerate}
\item
Let $X, Y$ be pointed compact metric spaces.
Then
$$d_{GH, p}\bigl(\cH(X), \cH(Y)\bigr) \le d_{GH, p}(X, Y).$$
\item
Let $X, Y$ be pointed compact metric spaces.
Then
$$d_{GH, p, c}\bigl(\cH(X), \cH(Y)\bigr) \le d_{GH, p, c}(X, Y).$$
\end{enumerate}
\end{cor}

\section{A remark on the Gromov--Hausdorff distance between dynamical systems, shifts on a torus}

\begin{dfn}
Let $X$ be a dynamical system.
We define a \emph{trajectory} (\emph{with the initial condition} $\hat{x}$) as the set $\{g\hat{x}\}_{g \in G}$.
\end{dfn}

Consider the torus $W = \bigl\{(\theta_1, \dots, \theta_{n}), \ \theta_{i} \mod 2\pi \bigr\}$ (with the Euclidean metric) and the \emph{shift on the torus}
$T_{\omega} \: W \to W, \ T_{\omega}(\theta) = \theta + \omega \mod 2\pi$.
Clearly, this map is an isometry and preserves the Lebesgue measure.
Let $W_{\omega} = (W, T_{\omega})$.

\begin{dfn}
The vector $\omega$ \emph{is called non-resonant} if for any $k \in \Z^{n} \setminus \{0\}$ we have $\langle k, \omega \rangle \notin 2\pi\Z$.
\end{dfn}

In other words, the vector $\omega = (\omega_{1}, \dots \omega_{n})$ is non-resonant if the vectors
$2\pi, \omega_{1}, \dots, \omega_{n}$ are linearly independent over $\Q$.
It is known (see, e.g.,~\cite{KornfeldSinaiFomin1980})
that the vector $\omega$ is non-resonant if and only if
the map $T_{\omega}$ is ergodic (which implies that almost all trajectories are everywhere dense).
Nevertheless, for any $\omega, \omega'$ the Gromov--Hausdorff distance $d_{GH, \cK}(W_{\omega}, W_{\omega'}) = 0$, which follows from a more general observation.

\begin{ass}
Let $X, Y$ be metric dynamical systems.
Suppose that for all $g \in G$ the maps $\phi_{X}(g, \bullet) \: X \to X$ and $\phi_{Y}(g, \bullet) \: Y \to Y$ are isometries.
Then the Gromov--Hausdorff distance $d^{G}_{GH, \cK}(X, Y)$ coincides with the Gromov--Hausdorff distance $d_{GH, \cK}(X, Y)$ between the corresponding metric spaces.
\end{ass}

\begin{proof}
Indeed, since the given maps are isometric, in the definition of the distortion and codistortion of maps (for dynamical systems) we will obtain the distortion of the map in the sense of metric spaces.
\end{proof}

In this regard, we introduce one more distance for dynamical systems.

To automatically obtain the triangle inequality for it, we will use the already constructed theory of $F$-spaces by associating a certain auxiliary space to each dynamical system.

To each dynamical system $X$ we associate the set $\Gamma_X = X \times G$. We define a single generalized pseudometric on it (that is, from the point of view of an $F$-space, the index set $I$ consists of a single element, $J = \emptyset$):
$$d_{\Gamma_{X}}\big((x_1, g_1), (x_2, g_2)\big) = d_X(g_1 x_1, g_2 x_2).$$
Note that an arbitrary map of dynamical systems $f \: X \to Y$ induces the map
$$\Gamma_{f} \: \Gamma_{X} \to \Gamma_{Y}, \ (x, g) \mapsto \bigl(f(x), g\bigr).$$

Based on the existing category $\cK \subset \cF$ we construct a category $\widehat{\cK}$ as follows.
By definition, the objects of this category $\widehat{\cK}$ are all $F$-spaces of the form $\Gamma_{X}$, where $X \in \cDGH$.
The morphisms are the maps:
$$\Mor_{\widehat{\cK}}(\Gamma_{X}, \Gamma_{Y}) =
\bigl\{ \Gamma_{f} \ | \ f \in \Mor_{\cK}(X, Y) \bigr\}.$$

Since by assumption the sets of all morphisms in $\cK$ are non-empty, the same is true for $\widehat{\cK}$.
Therefore, according to Proposition~\ref{ass:metric}, the categorical Gromov--Hausdorff distance $d_{GH,\widehat{\cK}}(\Gamma_X, \Gamma_Y)$ is a generalized pseudometric on $\Ob(\widehat{\cK})$. We will adopt this distance as the desired new distance between the dynamical systems themselves:
$$\widehat{d}^{G}_{GH, \cK}(X, Y) := d^{G}_{GH, \widehat{\cK}}(\Gamma_{X}, \Gamma_{Y}). $$

Let $f_1 \in \Mor_{\cK}(X, Y), \ f_2 \in \Mor_{\cK}(Y, X)$.
Let us compute the distortion and codistortion of the morphisms $\Gamma_{f_1}, \Gamma_{f_2}$. By the definition of distortion for $F$-spaces, we have:
$$ \widehat{\dis}^{G}(f_1) := \dis^{G}(\Gamma_{f_1}) = \sup_{g_1,g_2 \in G} \sup_{x_1,x_2 \in X} \Big| d_X(g_1 x_1, g_2 x_2) - d_Y\big(g_1 f_1(x_1), g_2 f_1(x_2)\big) \Big|. $$
$$ \widehat{\codis}^{G}(f_1, f_2) := \codis^{G}(\Gamma_{f_1}, \Gamma_{f_2}) = \sup_{g_1,g_2 \in G} \sup_{\substack{x \in X \\ y \in Y}} \Big| d_X\big(g_1 x, g_2 f_2(y)\big) - d_Y\big(g_1 f_1(x), g_2 y\big) \Big|. $$

Since the usual distortion and codistortion are obtained from this formula by setting $g_1 = g_2$ in the supremum, we have
$\widehat{\dis}^{G}(f) \ge \dis^{G}(f)$
and $\widehat{\codis}^{G}(f_1, f_2) \ge \codis^{G}(f_1, f_2)$,
which immediately implies the inequality
$\widehat{d}^{G}_{GH, \cK}(X, Y) \ge d^{G}_{GH,\cK}(X, Y)$.

Let us obtain a sufficient condition for reaching the upper bound in Proposition~\ref{prop:estimation} for the Gromov--Hausdorff distance $\widehat{d}^{G}_{GH}$.

\begin{thm}
Let $X, Y$ be dynamical systems.
Suppose there exists $\hat{x} \in X$ such that for all $\hat{y} \in Y$ the trajectory $\bigl\{(g\hat{x}, g\hat{y})\bigr\}_{g \in G}$ is everywhere dense in $X \times Y$.
Then
$$\widehat{d}^{G}_{GH, \cK}(X, Y) = \frac{1}{2} \max\bigl\{ \diam(X), \diam(Y) \bigr\}.$$
\end{thm}

\begin{proof}
We can assume that $X$ and $Y$ are bounded (the other cases are treated similarly).
Let $f_1 \in \Mor_{\cK}(X, Y), \ f_2 \in \Mor_{\cK}(Y, X)$.
Fix $\e > 0$.

Consider $x^{*}_1, x^{*}_2 \in X$ such that $d_{X}(x^{*}_1, x^{*}_2) > \diam(X) - \e$.
By the assumption, the trajectory $\bigl\{(g\hat{x}, gf_{1}(\hat{x}))\bigr\}_{g \in G}$ is everywhere dense in $X \times Y$.
We will assume that on $X \times Y$ we consider the generalized pseudometric of the $l^2$-product $X \times_{l^2} Y$.
Then there exists a point $g^{*} \in G$ such that
$$d_{X \times Y}\Bigl(\bigl(g^{*} \hat{x}, g^{*} f_{1}(\hat{x}) \bigr), \bigl(x^{*}_{2}, f_{1}(x^{*}_{1}) \bigr)\Bigr) < \e.$$
Then
$$d_{X}(g^{*} \hat{x}, x^{*}_{1}) \ge d_{X}(x^{*}_{1}, x^{*}_{2}) - d_{X}(g^{*} \hat{x}, x^{*}_{2}) \ge \diam(X) - 2\e
\text{ and }
d_{Y}(g^{*} f_{1}(\hat{x}), f_{1}(x^{*}_1)) < \e,$$
therefore
\begin{multline*}
\widehat{\dis}^{G}(f_1) = \sup_{g_1, g_2 \in G} \sup_{x_{1}, x_{2} \in X}
\Bigl| d_{X}(g_1 x_1, g_2 x_2) - d_{Y}\bigl(g_1 f_1(x_1), g_2 f_1(x_2)\bigr)\Bigr| \ge \\
\
\ge \sup_{g \in G} \sup_{x_{1}, x_{2} \in X}
\Bigl| d_{X}(g x_1, x_2) - d_{Y}\bigl(g f_1(x_1), f_1(x_2)\bigr)\Bigr| \ge        \\
\
\ge \Bigl| d_{X}(g^{*}\hat{x}, x^{*}_{1}) - d_{Y}\bigl(g^{*} f_{1}(\hat{x}), f_{1}(x^{*}_{1}) \bigr)\Bigr| \ge \diam(X) - 3\e.
\end{multline*}

Consider $y^{*}_{1}, y^{*}_{2} \in Y$ such that $d_{Y}(y^{*}_{1}, y^{*}_{2}) > \diam(Y) - \e$. Since the trajectory with the initial condition $\bigl(\hat{x}, f_{1}(\hat{x})\bigr)$ is everywhere dense in $X \times Y$, there exists a point $g^{*} \in G$ such that
$$d_{X \times Y}\Bigl(\bigl(g^{*} \hat{x}, g^{*} f_{1}(\hat{x}) \bigr), \bigl(f_{2}(y^{*}_{1}), y_{2} \bigr)\Bigr) < \e.$$
Then
$$d_{Y}\bigl(g^{*} f_{1}(\hat{x}), y^{*}_{1}\bigr) \ge
d_{Y}(y^{*}_{1}, y^{*}_{2}) - d_{Y}\bigl(g^{*} f_{1}(\hat{x}), y^{*}_{2}\bigr) \ge \diam(Y) - 2\e \text{ and }
d_{X}(g^{*}\hat{x}, f_{2}(y^{*}_{1})) < \e,$$
therefore
\begin{multline*}
\widehat{\codis}^{G}(f_{1}, f_{2}) =
\sup_{g_{1}, g_{2} \in G} \ \sup_{\substack{x \in X \\ y \in Y}} \
\Bigl| d_{X}\bigl(g_{1} x, g_{2} f_{2}(y)\bigr) -
d_{Y}\bigl(g_{1} f_{1}(x), g_{2} y\bigr) \Bigr| \ge \\
\
\ge \sup_{g \in G} \ \sup_{\substack{x \in X \\ y \in Y}} \
\Bigl| d_{X}\bigl(g x, f_{2}(y)\bigr) -
d_{Y}\bigl(g f_{1}(x), y\bigr) \Bigr| \ge
\Bigl| d_{X}\bigl(g^{*} \hat{x}, f_{2}(y^{*}_{1}) \bigr) -
d_{Y}\bigl(g^{*} f_{1}(\hat{x}), y^{*}_{1}\bigr)
\Bigr| \ge \\
\
\ge \diam(Y) - 3\e.
\end{multline*}
We have proved that for any $f_1 \in \Mor_{\cK}(X, Y), \ f_2 \in \Mor_{\cK}(Y, X)$ and any $\e > 0$ we have
\begin{multline*}
\max\bigl\{ \widehat{\dis}^{G}(f_1), \ \widehat{\dis}^{G}(f_2), \ \widehat{\codis}^{G}(f_1, f_2) \bigr\} \ge \\ \ge \max\bigl\{ \widehat{\dis}^{G}(f_1), \ \widehat{\codis}^{G}(f_1, f_2) \bigr\} \ge \max\bigl\{ \diam(X), \ \diam(Y) \bigr\} - 3\e.
\end{multline*}
By Proposition~\ref{prop:estimation} it follows that
$$\widehat{d}^{G}_{GH, \cK}(X, Y) = \frac{1}{2} \max\bigl\{ \diam(X), \ \diam(Y) \bigr\}.$$
\end{proof}

Let us return to our example: consider the shifts on the torus $W_{\omega}$ and $W_{\omega'}$, $\omega, \omega' \in \R^{n}$.
Consider the vector $(\omega, \omega') \in \R^{2n}$.
Note that $W_{(\omega, \omega')} = W_{\omega} \times_{l^{2}} W_{\omega'}$.
From the just proved theorem, we obtain the following corollary.
\begin{cor}
Suppose that the vector $(\omega, \omega')$ is non-resonant.
Then $$\widehat{d}^{G}_{GH, \cK}(W_{\omega}, W_{\omega'}) = \pi \sqrt{n}.$$
\end{cor}

\begin{proof}
Since the vector $(\omega, \omega')$ is non-resonant, all trajectories in $W_{(\omega, \omega')}$ are everywhere dense.
Then we obtain
$$\widehat{d}^{G}_{GH, \cK}(W_{\omega}, W_{\omega'}) =
\frac{1}{2} \max\bigl\{\diam(W_{\omega}), \ \diam(W_{\omega'}) \bigr\} = \pi \sqrt{n}.$$
\end{proof}

Thus, the new variant of the Gromov--Hausdorff distance $\widehat{d}^{G}_{GH,\cK}$ allows one to distinguish shifts on the torus.

 \end{document}